\documentclass[12pt]{article}
\usepackage{amsfonts}
\usepackage{epsfig}
\title{The Four Color Theorem meets Shapes of Polyhedra}
\author{Richard Evan Schwartz \thanks{\hskip 5 pt Supported by 
N.S.F. Grant DMS-2505281}}

\newtheorem{theorem}{Theorem}[section]

\newtheorem{lemma}[theorem]{Lemma}

\newtheorem{conjecture}[theorem]{Conjecture}

\def\startproof{{\bf {\medskip}{\noindent}Proof: }}

\def\endproof{$\spadesuit$  \newline}

\def\C{\mbox{\boldmath{$C$}}}%
\def\H{\mbox{\boldmath{$H$}}}%
\def\R{\mbox{\boldmath{$R$}}}%
\def\Z{\mbox{\boldmath{$Z$}}}%

\begin{document}

\maketitle

\begin{abstract}
  We consider solutions to the
  $4$-color problem for the vertices of 
  sphere triangulations with
  degree sequence $6,...,6,4,4,4,4,4,4$.
  We sort these solutions into combinatorial
  types and show that each generic type $\tau$ is
  parametrized by the set of integer lattice points inside a
  rational polyhedral convex cone ${\cal
    C\/}_{\tau}$ of real
 dimension is least $4$.
  There is an integral quadratic form $Q_{\tau}$
  on ${\cal C\/}_{\tau}$
  whose diagonal part, evaluated on a lattice point,
  is $3$ times the number of triangles in the corresponding
  triangulation.  We relate this structure to
  the octahedral stratum of Thurston's
  moduli space of flat cone structures on the
  sphere.
\end{abstract}

\section{Introduction}

\subsection{Context}

Here are two formulations of the famous Four Color Theorem.
\begin{enumerate}
 \item You can color the vertices of any sphere triangulation with $4$
   colors in such a way that
  no two adjacent vertices get the same color.
\item You can color the triangles of any sphere triangulation black
  and white, so that the numbers of black and white triangles
  incident to each vertex are congruent to each other mod $3$.
\end{enumerate}
Property (1)
defines a  piecewise affine map $\phi$ from the sphere to the
regular tetrahedron whose vertices have been given the same
$4$ colors. Color a triangle in the sphere triangulation white if and only
if $\phi$ is orientation preserving on that triangle.
This gives Property (2).  To go from Property (2) to Property (1)
you define the map $\phi$ using the orientation data and then
pull back the vertex coloring of the tetrahedron.  The mod $3$
condition guarantees that $\phi$ is well defined.

See [{\bf W\/}] for a lot of general information and
history about the Four Color Theorem.

Certain families of triangulations closely resemble moduli spaces,
and Property (2) above gives a way to study solutions to the
$4$ color problem on infinite collections of triangulations
simultaneously.  This paper  will explore this perspective in an
interesting special case, that of triangulations having
degree sequence $6,...,6,4,4,4,4,4,4$.
In some sense, this paper is a sequel to my paper [{\bf S1\/}], which
explores Property 2 relative to triangulations having
degree sequence $6,...,6,3,3,3$.  However, the aims of
this paper are different, and one can read this
paper independently from [{\bf S1\/}].

A sphere triangulation has {\it non-negative curvature\/} if
there are at most $6$ triangles incident to each vertex.
If you build such a triangulation out of equilateral
triangles you produce a metric space $\Sigma$ called a
{\it flat cone sphere\/}.  All but finitely many points of
$\Sigma$ are locally isometric to the plane, and the remaining
points, corresponding to vertices of degree less than $6$,
are locally isometric to Euclidean cones.  In the case
we are interested in here, $\Sigma$ has $6$ cone points,
all having cone-angle $4\pi/3$.   That is, $\Sigma$ is
intrinsically an equiangular octahedron.

In his paper {\it Shapes of Polyhedra\/} [{\bf T\/}], Bill
Thurston organized the triangulations of non-negative
curvature along these lines.  Thurston fixes the list of $k$ cone points
and corresponding cone angles, and then considers flat
cone spheres with these cone points.
The result is a moduli space
$\cal M$ of flat cone spheres having prescribed (and,
for ease of discussion, labeled) cone points.
The space $\cal M$ is locally modeled on $\C^{1,k-3}$,
a complex vector space with a Hermitian form of
signature $(1,k-3)$.  Moreover, $\cal M$ has a
system of coordinate charts in which the
points with coordinates in the Eisenstein lattice
correspond to triangulations.
When one considers these
flat cone spheres modulo complex scaling,
which amounts to projectivizing the moduli space, the
resulting space is a complex hyperbolic manifold
of complex dimension $k-3$.

From now on we set
\begin{equation}
  {\cal M\/}={\cal M\/}(4,4,4,4,4,4).
\end{equation}
The space $\cal M$ is locally modeled on
$\C^{1,3}$.  The projectivized space
$P\cal M$ is 
an open
complex hyperbolic $3$-manifold whose metric
completion is the quotient of $\C\H^3$
by a certain lattice known as a
{\it Mostow-Deligne lattice\/}.

\subsection{Nice Colorings}

We say that a {\it nice polygon\/} is any convex polygon whose
interior angles all have the form $k\pi/3$ for $k\in\{1,2\}$.
A nice polygon has at most $6$ sides.  There are $5$ kinds:
triangles, parallelograms, trapezoids, pentagons and hexagons.
We call a nice polygon {\it unit triangulable\/} if it is tiled by
equilateral triangles having unit integer side lengths.
For instance, the regular hexagon, chosen to have
unit side lengths, is unit triangulable.

We say that a {\it nice coloring\/} of a flat cone sphere $\Sigma$ in
${\cal M\/}$ is a partition of $\Sigma$ into
nice polygons, alternately colored black and white, such that
\begin{enumerate}
\item There are $4$ polygons around each regular vertex.
\item There are $2$ polygons around each cone vertex.
\end{enumerate}
Here the {\it regular vertices\/} are the vertices of the partition
that are contained in disk neighborhoods.
We call the nice coloring {\it unit triangulable\/} if the polygons
in the partition are all unit triangulable.  In this case, the
coloring gives a solution of the $4$ coloring problem for
the corresponding sphere triangulation.

Figure 1.1 shows an example of a nice coloring.
In this example, the edges of the
polygon glue together to make the flat cone sphere.
Each boundary edge is paired isometrically to another
boundary edge in a direction-reversing way.  For instance,
the rightmost two edges are glued together to create one
of the cone points.  Figure 2.1 shows how the gluing works
in this example.

\begin{figure}[htbp]
\centering
\includegraphics[height=2.5in]{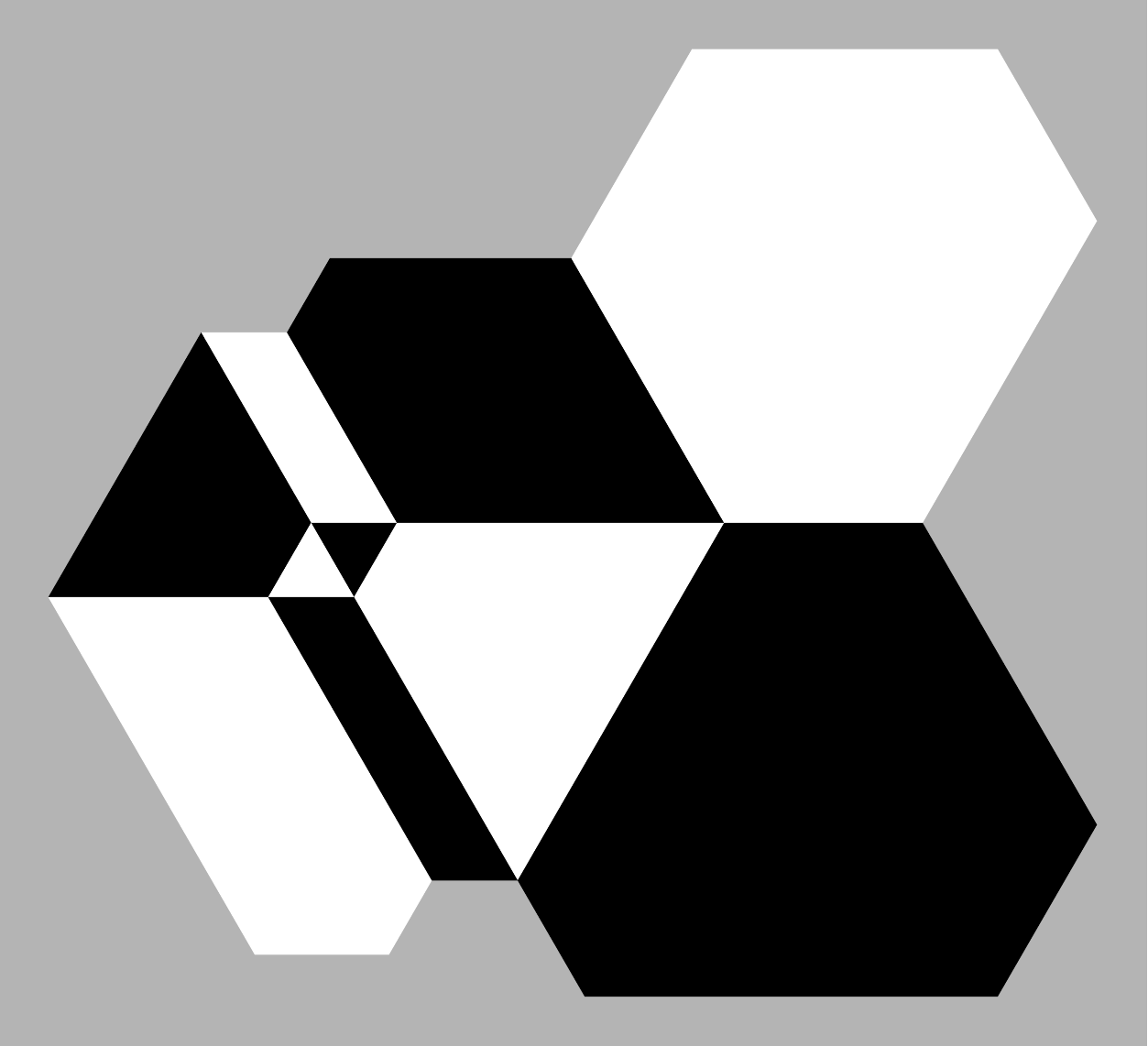}
\caption{A nice coloring of a flat cone sphere.}
\end{figure}

Figures 1.2-1.4 show the example from Figure 1.1 and three other
examples having the same combinatorics and gluing pattern.

\begin{figure}[htbp]
\centering
\includegraphics[height=2.5in]{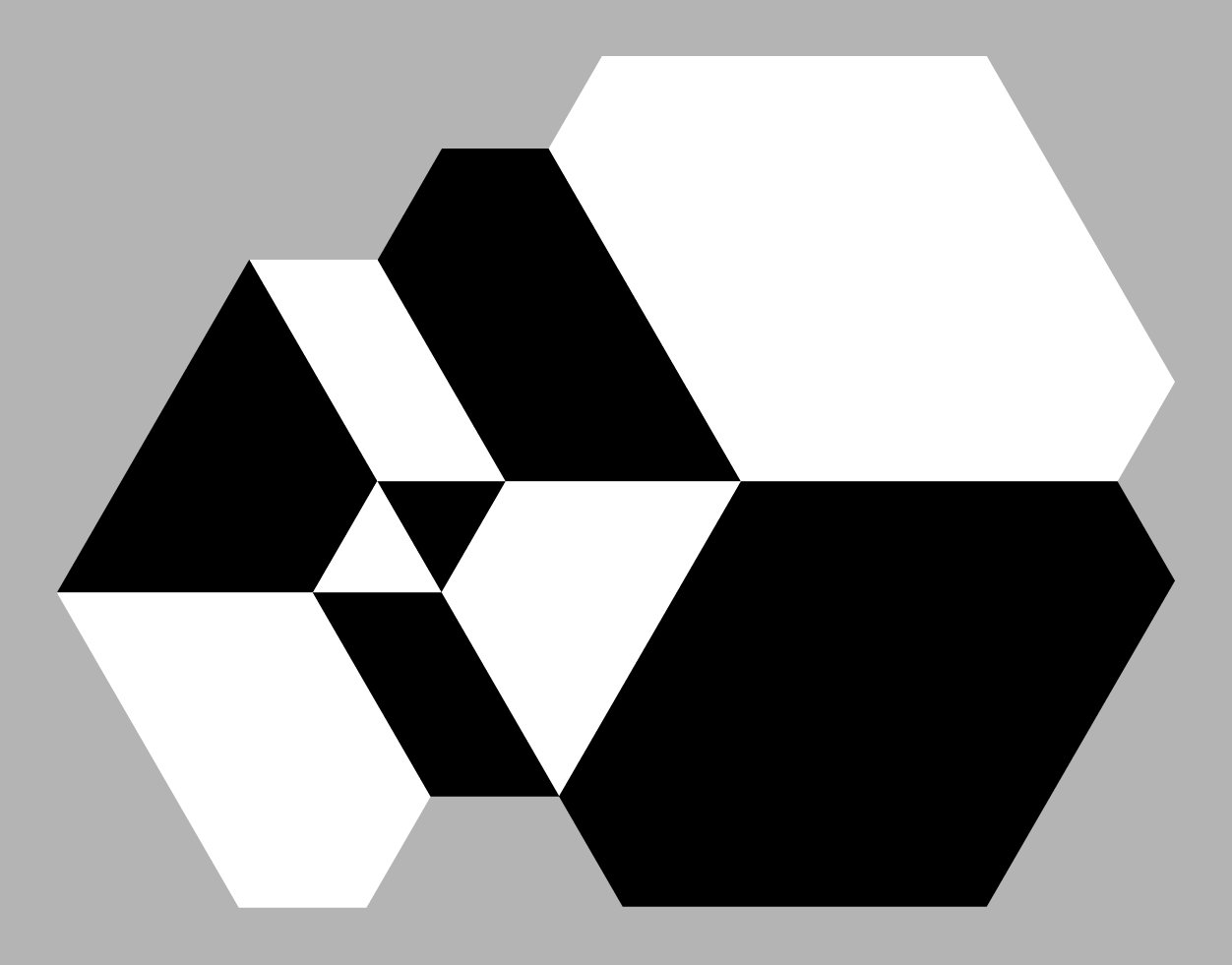}
\caption{Same combinatorics, different geometry}
\end{figure}

\begin{figure}[htbp]
\centering
\includegraphics[height=2.5in]{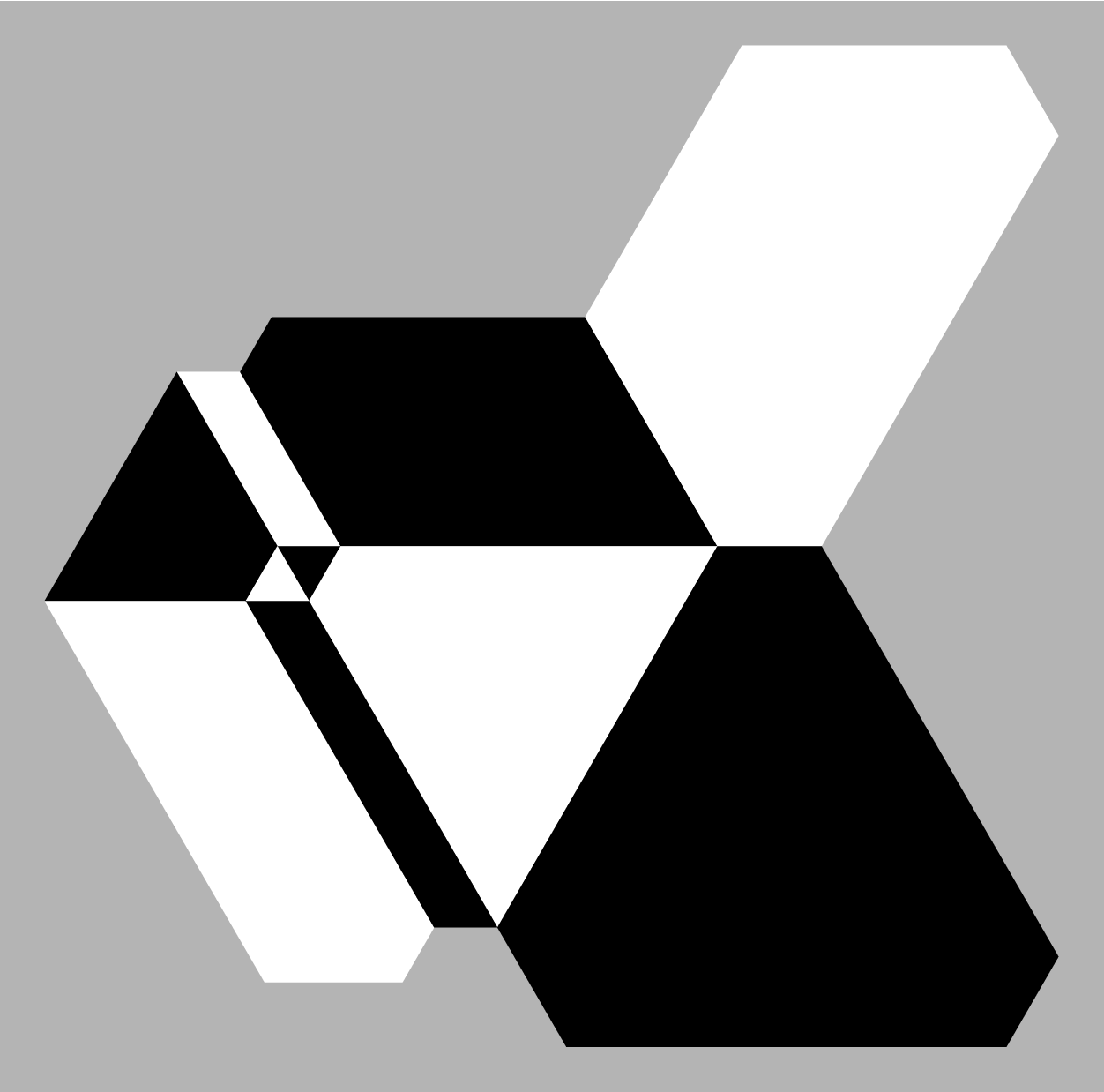}
\caption{Same combinatorics, different geometry}
\end{figure}

\begin{figure}[htbp]
\centering
\includegraphics[height=2.5in]{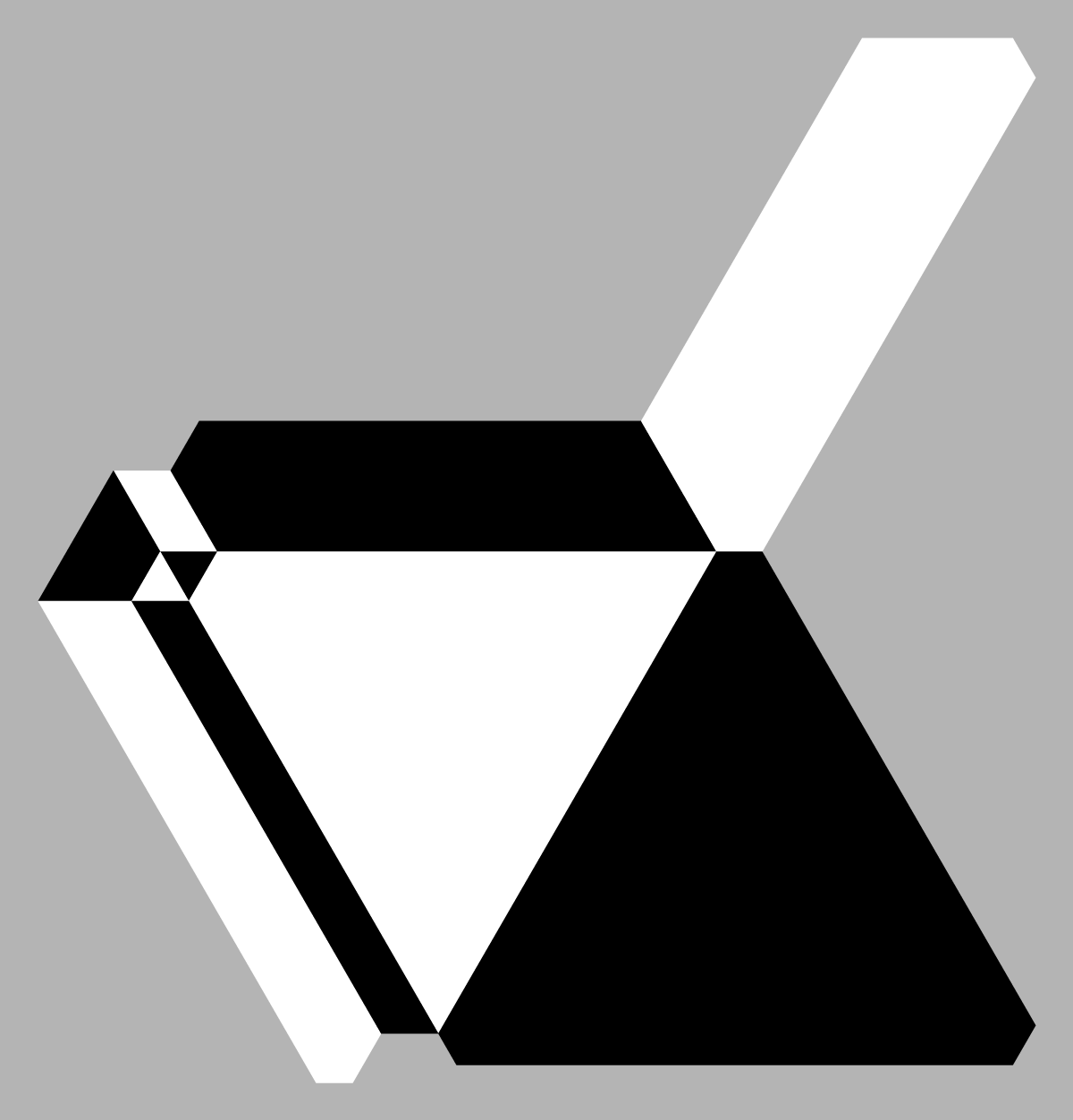}
\caption{Same combinatorics, different geometry}
\end{figure}

The examples shown in Figures 1.2-1.4 come from a $4$-parameter
family of examples, naturally parametrized by a $4$-dimensional
convex cone ${\cal C\/} \subset \R^E$, where $E$ is the number of
edges.
The cone $\cal C$ is a rational polyhedral convex cone.
The $4$-dimensional linear subspace extending
${\cal C\/}$ intersects $\Z^E$ in a co-compact lattice.
The points of ${\cal C\/} \cap \Z^E$ give rise to
unit triangulable nice colorings, all having the
same combinatorial structure.
In this way,
the same combinatorial pattern describes an infinite
family of solutions to the $4$ color problem.
The projectivized space $P\cal C$ is a
$3$-dimensional convex polytope, and it has
a dense set of points corresponding to
triangulable nice colorings.
\newline
\newline
{\bf Remark:\/}
The precise local conditions we gave for a nice coloring
do not cover all possible colorings which correspond to
solutions of the $4$-color problem.  The other kinds of
colorings by nice polygons (with more general
kinds of local pictures) also come in families, but
they are lower dimensional.   This is why we say that
the $4$-color solutions corresponding to nice
colorings are the generic kind.   For the sake of
simplicity, we do not treat the other kinds of
$4$-coloring solutions even though we could
use similar methods.

\subsection{Main Results}

Our first result generalizes the discussion surrounding
the example from Figure 1.2.

\begin{theorem}
  \label{one}
  Let $C$ denote any nice coloring of a flat cone
  octahedron. Let $E$ denote the number of edges
  in the coloring.  The set of nice colorings combinatorially
  equivalent to $C$ is parametrized by a
  rational  polyhedral convex cone ${\cal C\/} \subset \R^E$.
  The subset of triangulable nice colorings combinatorially
  equivalent to $C$ is
   parametrized by  ${\cal C\/} \cap \Z^E$.
\end{theorem}
In \S \ref{realize},
we will give a concrete recipe for
constructing $\cal C$ given the combinatorics
of $C$.   The construction can be implemented
by computer, and in fact my computer program
implements it.

The nice colorings seem to be quite abundant,
and some are quite beautiful.  In \S 2 I will show
a lot of examples.   The two most interesting
examples I found are each indexed by
the set of finite integer  sequences $a_1,...,a_k$
where $a_i \geq 2$ for all $i=1,...,k$, and
$k \geq 1$ can be as large as we like.
I have not tried to prove rigorously that these families
exist, but the numerical evidence
is overwhelming.   The existence of these
families would show that the number of
  combinatorially distinct nice colorings
  with $2n$ polygons grows at least
  exponentially in $n$.
  \newline
  
  We call $\cal C$ a {\it cone of colorings\/}.
  We set ${\cal C\/}_{\Z}={\cal C\/} \cap \Z^{E}$.
Our next result relates $\cal C$ to Thurston's
moduli space $\cal M$.  Just as Thurston has
a natural Hermitian form on his moduli spaces,
we will have a natural quadratic form.
Our quadratic form is very reminiscent of
the one in [{\bf BG\/}].

In all the examples I have computed for
this paper, $\cal C$ is
$4$ dimensional.   Dia Taha and Kai Fu gave a
simple example showing that the dimension
can exceed $4$.  One can soup up their example
to show that the dimension of $\cal C$ can
be arbitrarily large.  These examples have a
special kind of degeneracy.

\begin{theorem}
  \label{two}
  Suppose that the cone $\cal C$ of colorings
  is $4$-dimensional.
There is an integral quadratic form $Q$
  on $\cal C$ and a locally affine map
  $A: {\cal C\/} \to {\cal M\/}$.   For each
  $V \in {\cal C\/}_{\Z}$, the expression $Q(V,V)$ computes
  $3$ times the number of triangles in the
  triangulation associated to $V$.
  If $Q$ is a non-degenerate form then
  $A$ is locally injective.
    If $Q$ has signature $(1,3)$ then
  the projectivized map $PA$ from $P{\cal C\/}$
  to $P{\cal M\/}$ is locally
  injective. 
\end{theorem}

The quadratic form $Q$ turns out to be the
pullback of the real part of Thurston's
Hermitian [{\bf T\/}] form when it is suitably
normalized.   We will give a direct formula for
$Q$ in \S \ref{QF}.  The formula is completely
combinatorial; it does not depend on the
geometry of any particular nice coloring
associated to $\cal C$.
Theorem \ref{two} would be more
powerful with the following conjecture.

\begin{conjecture}
  \label{signature}
  When the cone of colorings is $4$-dimensional,
  the quadratic form $Q$ always has
  signature $(1,3)$.
  In particular, $Q$ is always non-degenerate.
\end{conjecture}

The conjecture has some nice payoffs.
First, if $Q$ has signature $(1,3)$ then
$Q$ imparts a real hyperbolic structure
to $P\cal C$. Second,
the conjecture combines with Theorem \ref{two} to show that
the polyhedron $P\cal C$ always maps in a locally
injective way into $P\cal M$.  
Conjecture \ref{signature}
lines up with the situation in
[{\bf BG\/}] and [{\bf T\/}], but I don't think that
the conjecture follows from anything in those
papers.
In any given case I can compute $Q$
in a basis for the subspace spanned by $\cal C$.
In all the many cases I have checked, the signature
is $(1,3)$.

Here is a related conjecture.
\begin{conjecture}
  When $\cal C$ is $4$-dimensional,
  the map $PA: P{\cal C\/} \to P{\cal M\/}$ is (globally)
  injective.
\end{conjecture}

This paper is still preliminary.  In any case, there
is plenty more to say.  Here are some additional
topics I might add, either to this paper or to a
sequel.
\begin{itemize}
\item Explicit calculations of the quadratic form on some
  infinite families of examples.  In some of the families,
  and with respect to a suitably chosen family of bases,
  the sequence of quadratic forms converges to a
  limiting form, and the limit (also) has signature $(1,3)$.
\item A discussion of the other combinatorial
  kinds of colorings of a flat cone octahedron
  by nice polygons.
\item A discussion of how the various convex
  cones from Theorem \ref{one} fit together, and
  relatedly how the nice coloring degenerates
  when we move to a face of the corresponding
  convex cone.
\end{itemize}
Also, I still need to make the software for
this paper public.   Right now, it lacks a
lot of documentation.  I hope to fix this
within a few weeks.

\subsection{Overview}

This paper is organized as follows.
In the beginning of
\S 2 I will define the basic combinatorial
objects of study, which is essentially the
multigraph dual to the nice coloring.
Following this, in \S \ref{simpleI}, I
will exhibit the infinite family of nice colorings.
I essentially give a proof, using a grafting
idea, but I stop short of giving all the details.
After this, I will present a gallery of
other examples of nice colorings.
In \S 3 I will prove Theorem \ref{one}.
In \S 4 I will prove Theorem \ref{two}.

\subsection{Acknowledgements}

I would like to thank Dia Taha, Kai Fu, and Emma Ruter
for various conversations about this project.  I would
also like to acknowledge that I used ChatGPT 5.6 Sol to check
the most recent version of the paper.  ChatGPT helped
me understand an issue with Lemma 3.3 in the
previous version, and it also suggested a quicker proof
of Lemma 3.3.

\newpage

\section{A Gallery of Examples}

\subsection{The Dual Multigraph}
\label{dualM}

In this section we explain how to encode our nice
colorings as a multigraph.  For us, an {\it octahedron\/} is
 always a flat cone sphere with $6$ cone points, all
 having cone angle $4\pi/3$.  Also, our {\it graphs\/} are
 allowed to have multiple edges between vertices.   Sometimes
 I will say {\it graph\/} and sometimes {\it multigraph\/}.

Let $P$ be an octahedron with a nice coloring.
We form a {\it dual multigraph\/}  $G$ as follows.
There is one vertex of $G$ for each polygon in
the triangulation.  We join two vertices by one edge
for each of the edges they have in common.
We draw this graph on $P$ by choosing the
vertices to be the barycenters of the faces and
by drawing each (blue) edge of the graph so that it
only crosses the corresponding edge of the
coloring of $P$.  Figure 2.1
shows an example.  Here, the
parts of the blue arcs which stray outside the
black and white polygon would disappear when
the polygon is glued up to make the octahedron.

\begin{center}
\resizebox{!}{4in}{\includegraphics{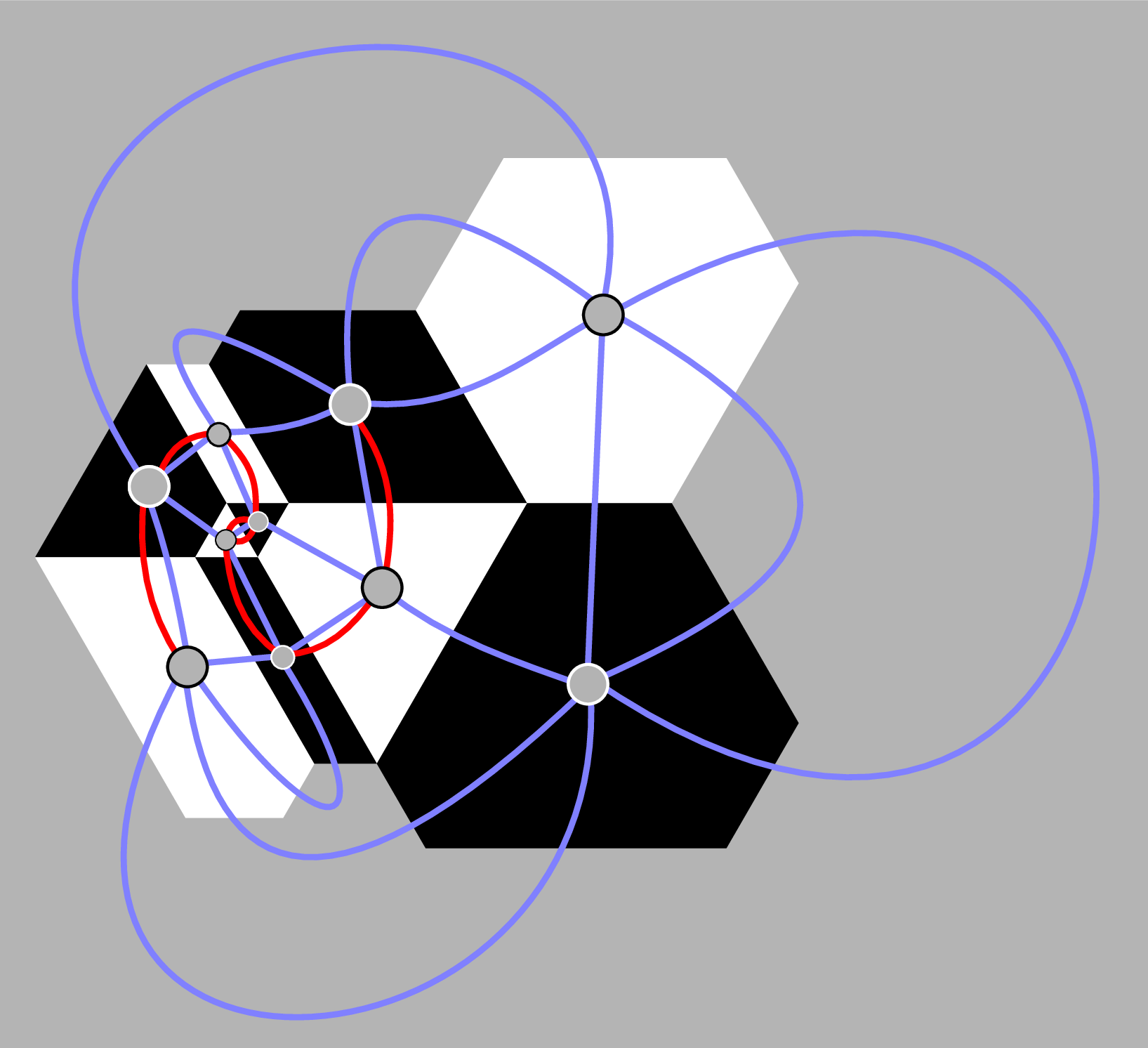}}
\newline
{\bf Figure 2.1:\/} The dual multigraph
\end{center}

This construction gives $G$ the
structure of an embedded planar multigraph.
The regions of $G$ bounded by bigons are in bijection
with the cone points of $P$.  So, there are $6$ bigon
regions in total. The remaining regions are quadrilaterals,
and they are naturally in bijection with the ordinary
vertices of $P$.  Figure 2.1 shows the example
that corresponds to Figure 1.1.

Each quadrilateral face has a distinguished edge $e$, the
one corresponding to the two acute angles around the
vertex corresponding to the face.  We add a
red edge to $G$ just inside the face and essentially
parallel to $e$. Figure 1.1 shows
these edges drawn in red.   When we add in all these
red edges we get what we call the {\it enhanced multigraph\/},
and we denote it by $\widehat G$.

\begin{lemma}
  Each vertex in $\widehat G$ has $6$ edges
  incident to it.
\end{lemma}

\startproof
Each vertex $v$ in the multigraph $G$,
corresponding to a nice $k$-gon $P$, has
$k$ blue edges incident to it,
one per side of $P$.  We want to see that
there are $6-k$ additional red edges incident to $v$.
We call a vertex $w$ of $P$ {\it acute\/} if its interior
angle at $w$ is acute.

\begin{center}
\resizebox{!}{1.4in}{\includegraphics{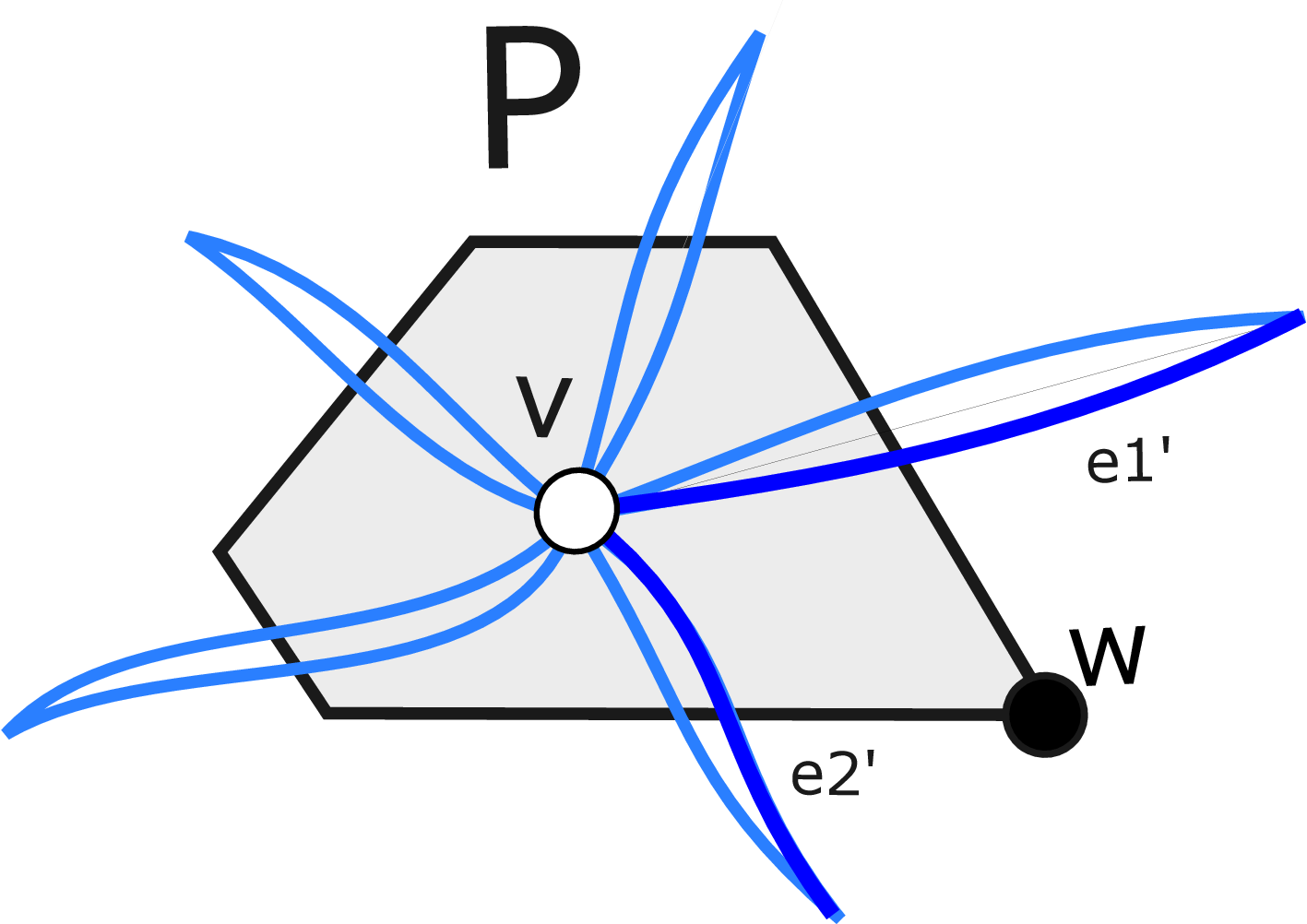}}
\newline
{\bf Figure 2.2:\/} The dual multigraph
\end{center}

We split each blue edge $e$ incident to $v$ into two infinitesimally
close {\it half-edges\/}.  Figure 2.2 shows this
schematically. Thus, there are $2k$ half-edges
emanating from $v$ and they are cyclically ordered.
Each half edge $e'$ defines a unique vertex $w$ of
$P$.  One can connect $e'$ to $w$ by an arc that
crosses no other blue edges of $G$.
We call two half-edges {\it partners\/} if they are
associated to the same vertex $w$ of $P$.  We call
these partners {\it acute\/} if the interior
angle at $w$  is acute.  Note that
there are exactly $6-k$ acute partners.

Given any pair $(e_1',e_2')$ of acute partners,
let $w$ be the associated vertex.
There are three other nice polygons
incident to $w$, and exactly one of
them has $w$ as an acute vertex.
Thus, exactly one of our half edges
in the set $\{e_1',e_2'\}$ has a parallel
red edge.  This gives us one red edge
for each pair of acute partners.  Hence there
are $6-k$ red edges incident to $v$.
\endproof

\subsection{A Simple Infinite Family}
\label{simpleI}

Figure 2.3 shows the beginning of a
family of cell divisions of
the sphere into quadrilaterals.  Our convention
is that the outside of the big square is also
part of the cell division.

\begin{center}
\resizebox{!}{1.25in}{\includegraphics{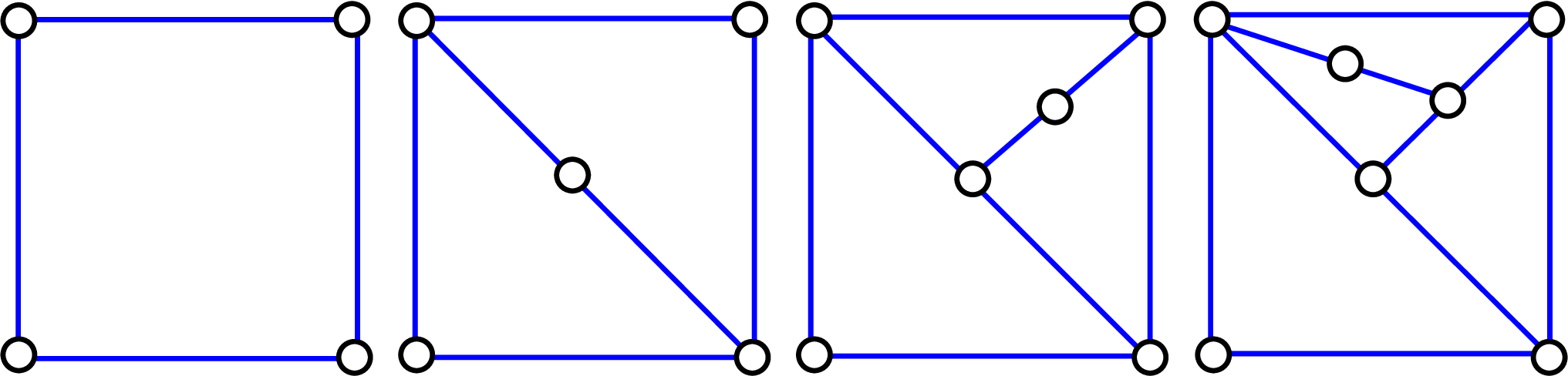}}
\newline
{\bf Figure 2.3:\/}  A family of cell divisions of the sphere
\end{center}
Thus, these
cell divisions respectively have $2,3,4,5$ faces.
Likewise they have $4,5,6,7$ vertices.
We are interested in the ones which have
an even number of vertices.   These are bipartite.
Based on the even members of this family, we
create red/blue multigraphs
$\widehat G_6, \widehat G_8, \widehat G_{10},...$
as shown in Figure 2.4.

\begin{center}
\resizebox{!}{1.5in}{\includegraphics{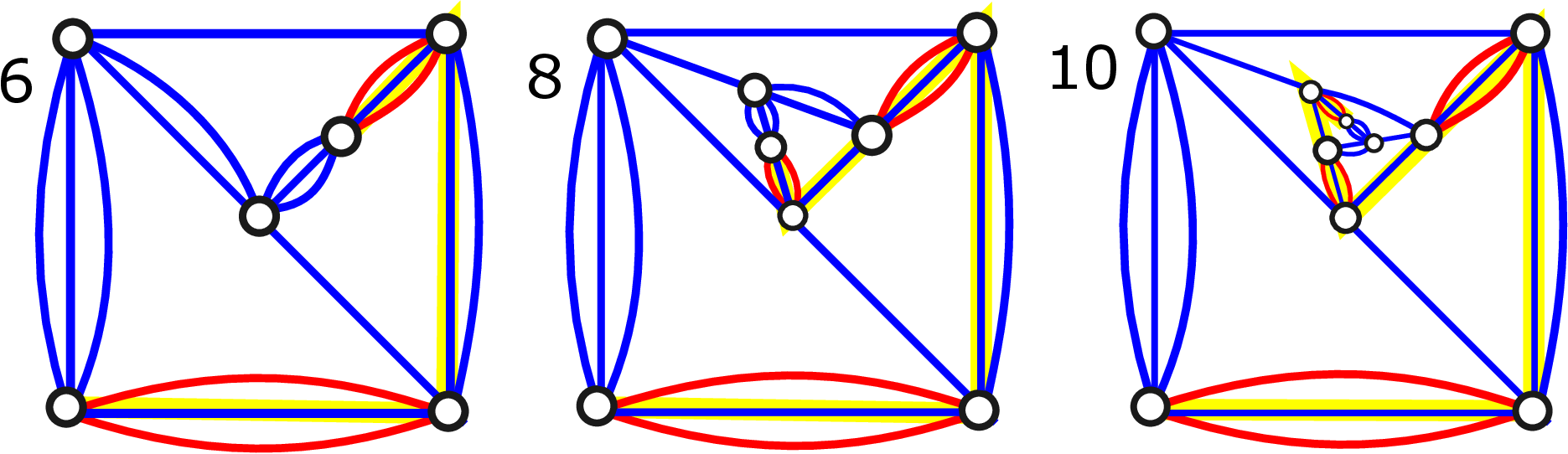}}
\newline
{\bf Figure 2.4:\/}  The graphs $\widehat G_6,\widehat G_8,\widehat G_{10}$
\end{center}

The continuation of the pattern may not be entirely clear, so we
say more about this.  There is a spiral path in $G_{2k}$ which
starts at the bottom left vertex and moves around until it reaches the
last placed vertex.  We have highlighted this path in yellow
in each case.  We put double red edges along this path, starting
with the bottom edge and then continuing in a way that
skips over every other edge.  We then add $6$ more blue
edges to make all the degrees $6$.

Figure 2.5 shows two nice colorings associated to
$\widehat G_{16}$. For the one on the left, we
have indicated the gluing pattern
both by letters and by some red zig-zags.
The one on the right has the same gluing
pattern.  We included both pictures to indicate
somewhat how the shapes change with the
parameter in the cone $\cal C$ guaranteed
by Theorem \ref{one}.

\begin{center}
\resizebox{!}{3in}{\includegraphics{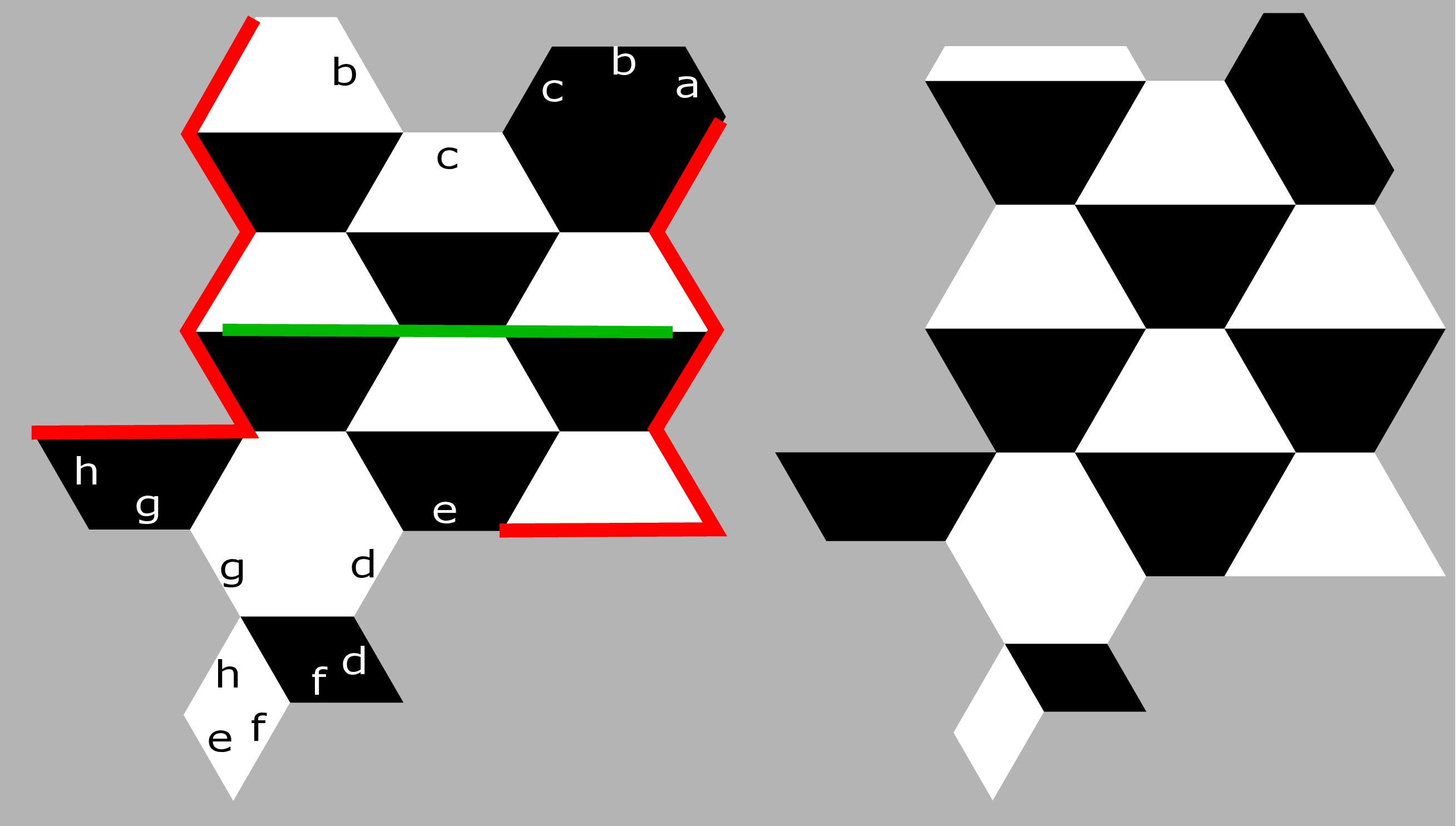}}
\newline
{\bf Figure 2.5:\/} A nice coloring associated to $\widehat G_{16}$.
\end{center}

Notice the red zigzags on either side of the figure on the left.
These are glued together by a translation.  This creates
a flat cylinder inside the octahedron.   We can lengthen
this cylinder by grafting in a $2 \times 3$ layer
of trapezoids.  We slit open the coloring along the
green segment, then make the graft using
trapezoids which have precisely the same shape as
the ones above and below the green slit.
We then modify the gluings so that longer
red zigzags are identified by a translation.
This produces a nice coloring with a longer cylinder.
This bigger nice coloring corresponds to
$\widehat G_{22}$.    This gives us a
provably infinite family of nice colorings.

The rest of this chapter presents a number of examples
without proof.  In each case, we will show one nice
coloring and its associated dual multigraph.
We include this material to illustrate the richness
of the construction.

\subsection{A Checkerboard Family}

Figure 2.6 shows one example from
an infinite family.  The gluing around the
boundary is very much like the folding of a wallet
or a taco.  There is an orientation-reversing
isometry of the boundary which fixes the two
vertices marked by the two dots, and then the
rest of the gluing is implemented by this
isometry.  A grafting operation like the one
presented in the previous section would
construct infinitely many examples from
this one.

\begin{center}
\resizebox{!}{5in}{\includegraphics{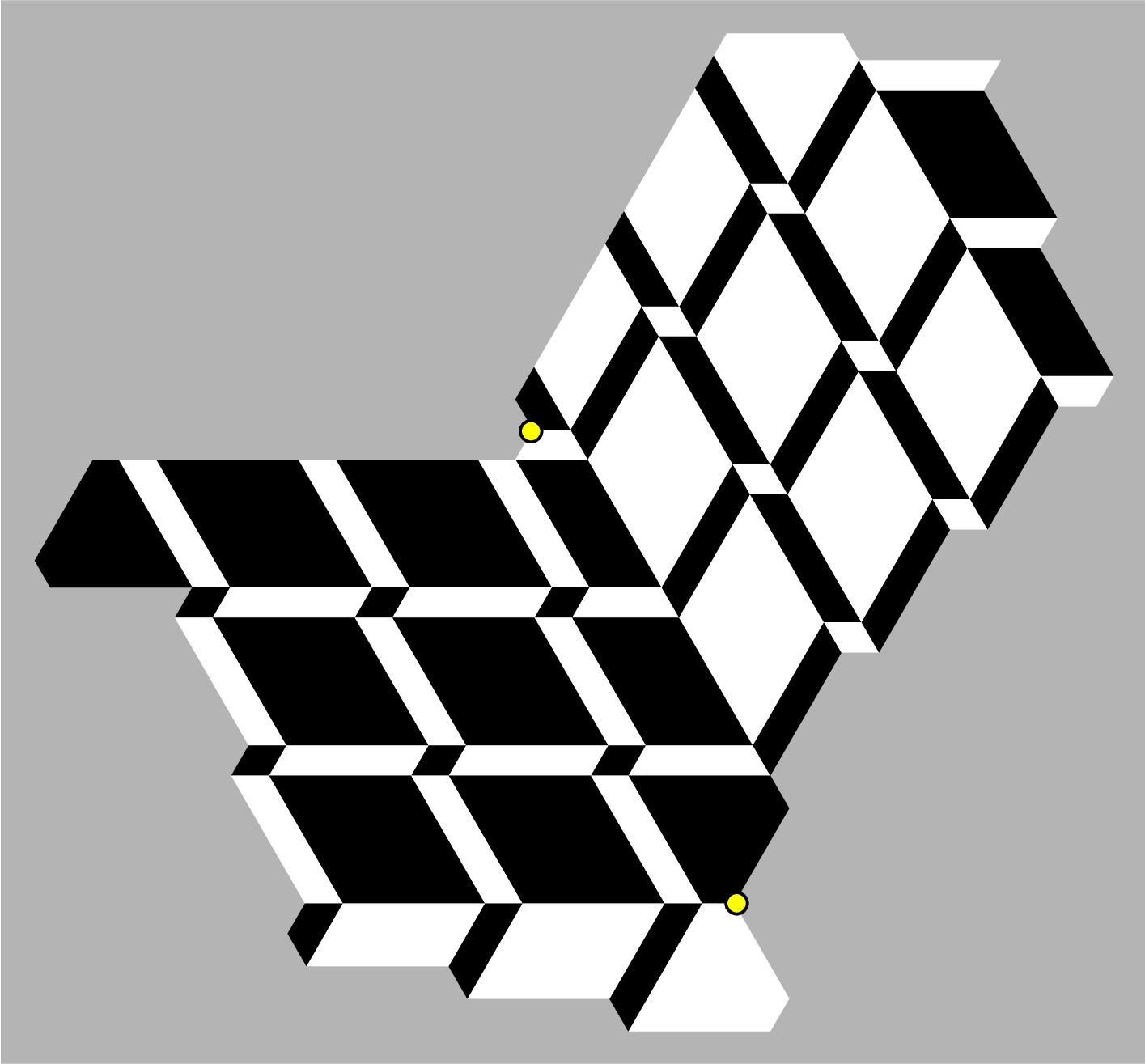}}
\newline
{\bf Figure 2.6:\/} A checkerboard themed example
\end{center}

Figure 2.7 shows the multigraph on which this example
is based.   The graph is made by gluing two identically
decorated disks together with an orientation reversing
isometry, in such a way that the arrowed edges
on the upper left of each side are identified.  A similar pattern like this
exists for any grid with an even number of squares,
though I did not explicitly check that all such graphs
actually correspond to nice colorings.

\begin{center}
\resizebox{!}{2.4in}{\includegraphics{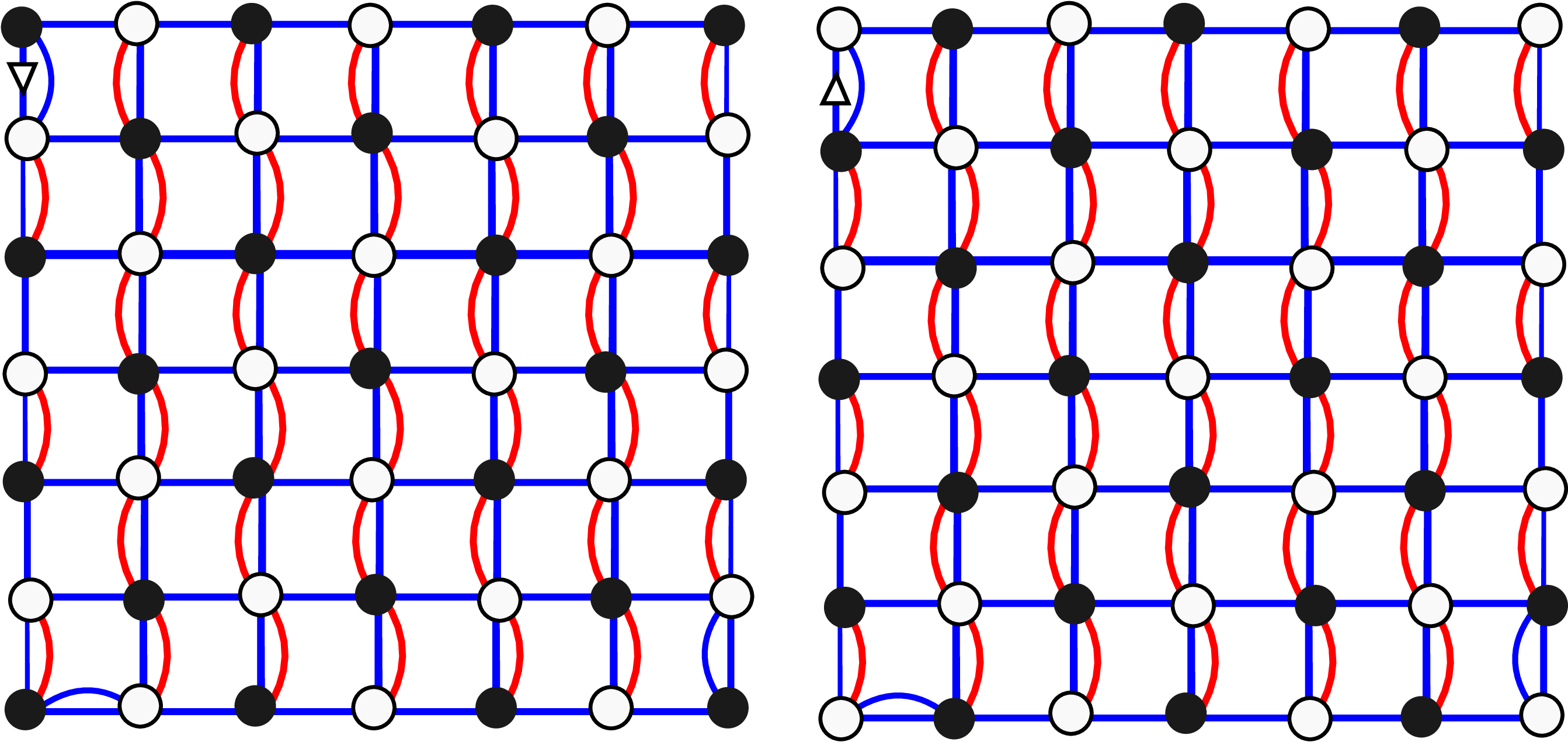}}
\newline
{\bf Figure 2.7:\/} The dual multigraph associated to Figure 2.6.
\end{center}

\subsection{A Lamination Themed Family}

Figure 2.8 shows two examples from what is presumably
an infinite family.  Unlike the previous two families, I
don't see an easy way to prove that all the members
in the family exist.  These examples (to me) give a
hint of geodesic laminations.

\begin{center}
\resizebox{!}{2.7in}{\includegraphics{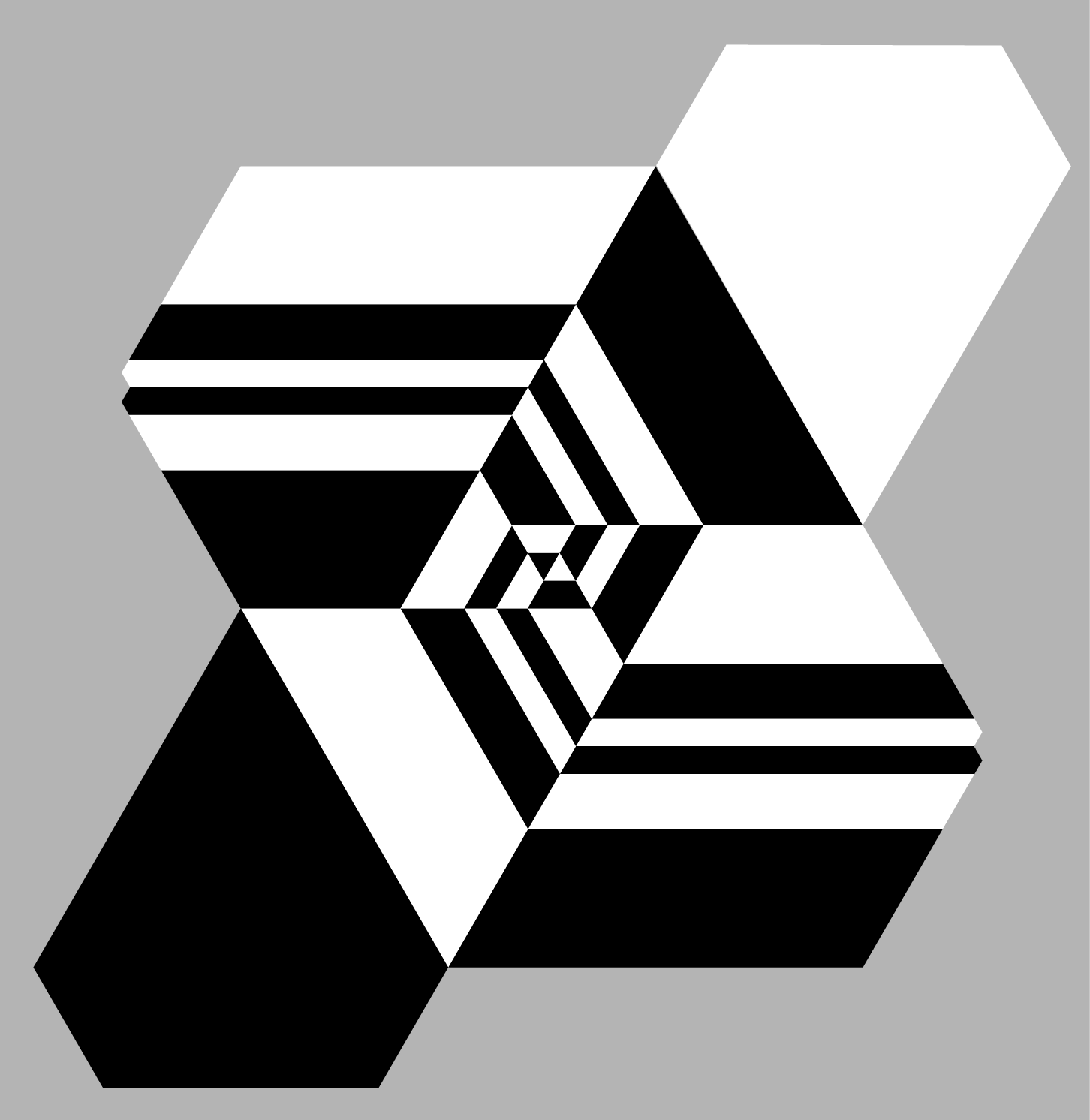}}
\resizebox{!}{2.7in}{\includegraphics{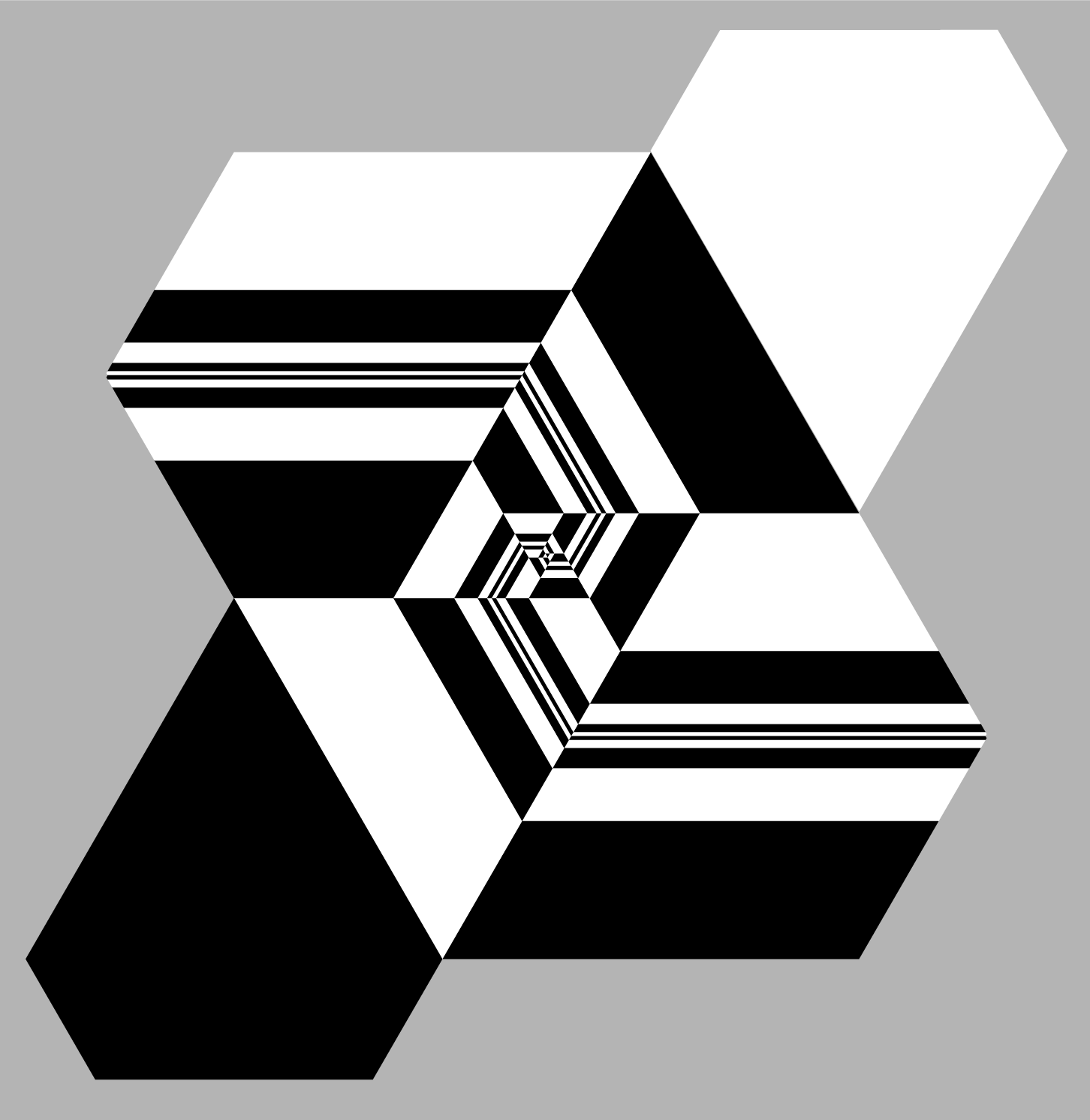}}
\newline
{\bf Figure 2.8:\/} Nice colorings which hint at geodesic laminations
\end{center}

Figure 2.9 shows the dual multigraph associated
to the left side of Figure 2.8. 
 In this picture we glue together
two squares in the pattern shown to make the sphere.
The lettering indicates the gluing pattern.
(I made the identifications more explicit in
this example so that the reader could trace
out the paths made by the red edges.)
For what it is worth,
I have used several shades of red to indicate the $3$
connected components of the red subgraph.

\begin{center}
\resizebox{!}{2.5in}{\includegraphics{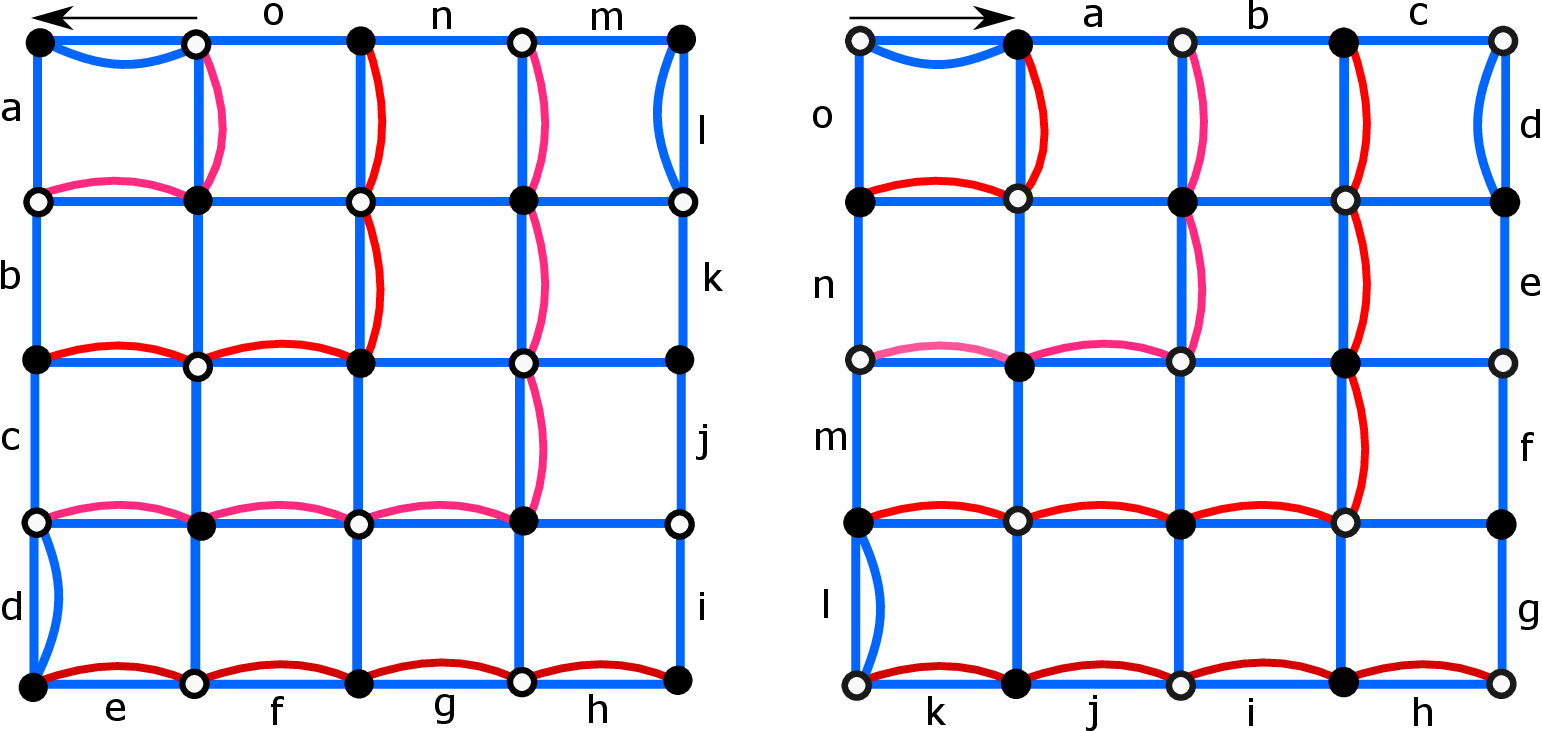}}
\newline
{\bf Figure 2.9:\/} The multigraph corresponding to the left side of
Figure 2.8
\end{center}

\subsection{Two Staircase Families}

I will illustrate the two big families with examples.
The examples are based on the sequence $2,3,2,4$.
These numbers control the successive lengths of a
path of squares which always moves upwards and to the right.
Figure 2.10 shows the dual graph for the corresponding
member of the first family.

\begin{center}
\resizebox{!}{1.4in}{\includegraphics{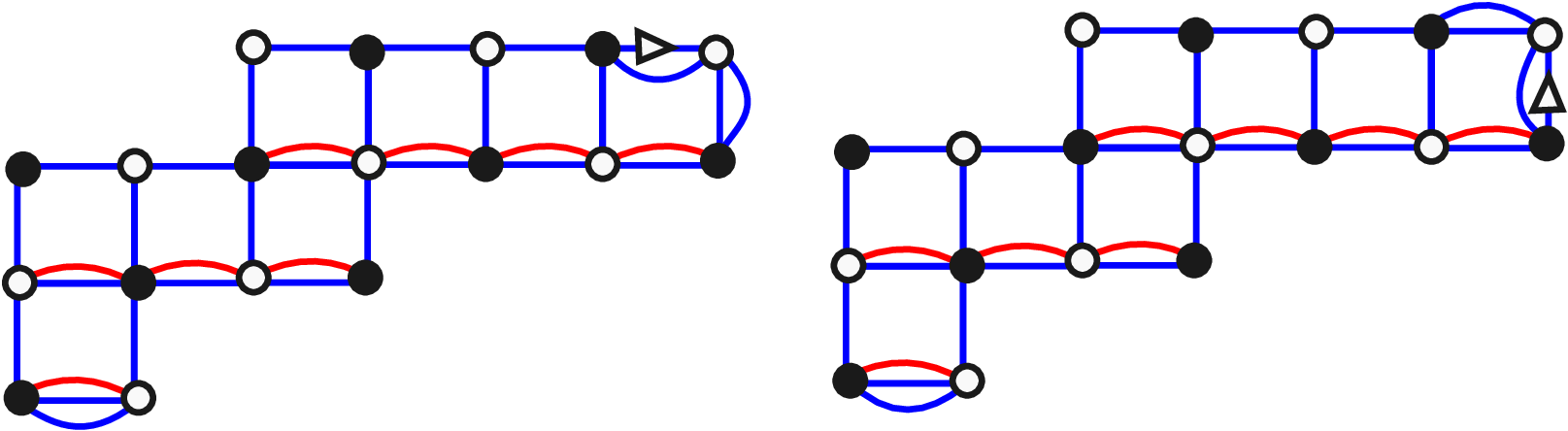}}
\newline
{\bf Figure 2.10:\/} Member $2324$ of the first staircase family
\end{center}

Figure 2.11 shows the corresponding nice coloring.
One thing I should say about my program is that it
does not necessarily lay out the most efficient
gluing diagram.  After computing the shapes,
according to the method described in the next chapter,
the program lays the shapes down in a kind of
random order, snapping them on to a growing
union of shapes. 

\begin{center}
\resizebox{!}{3.9in}{\includegraphics{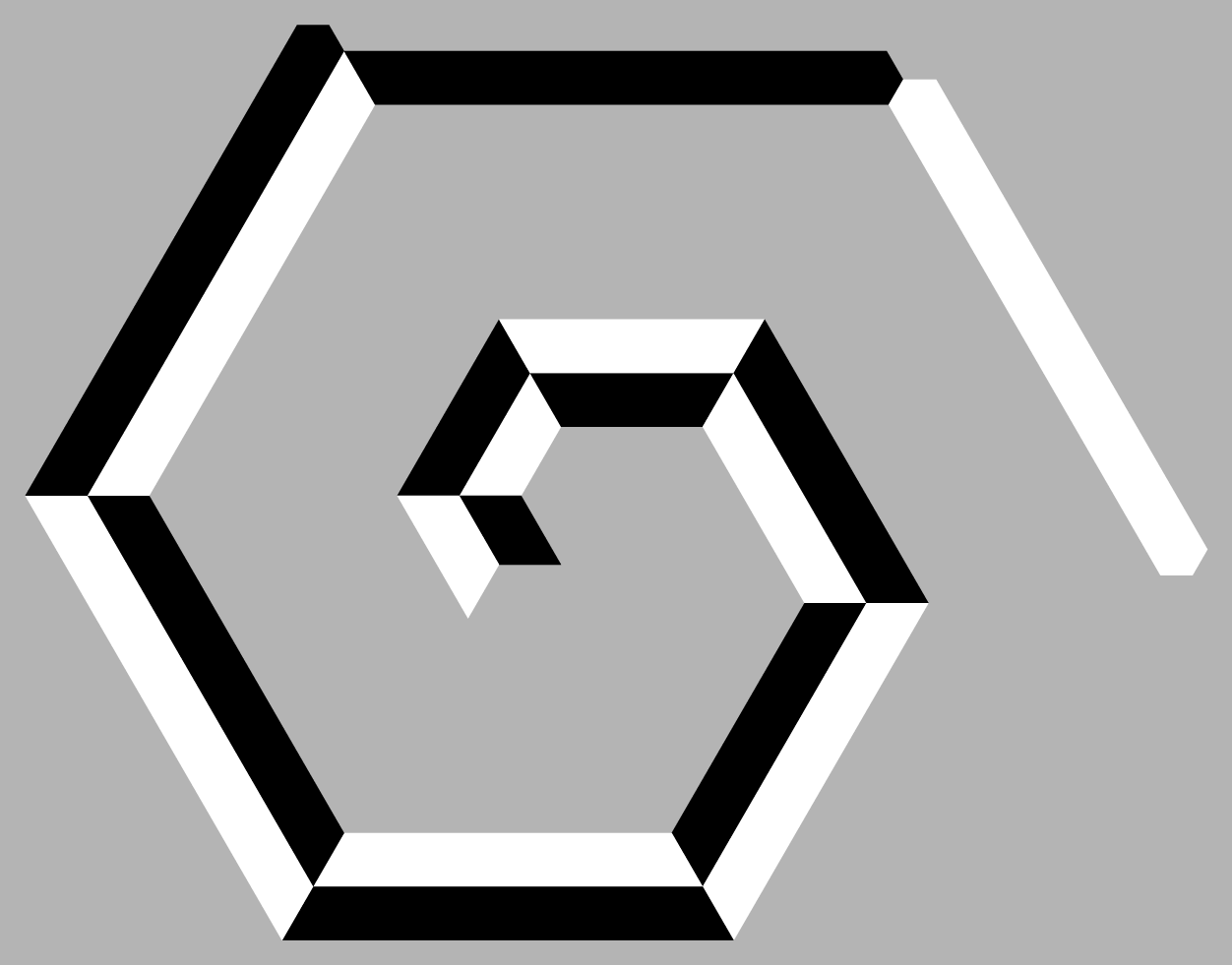}}
\newline
{\bf Figure 2.11:\/} The coloring associated to Figure 2.10.
\end{center}
I am not showing the explicit folding pattern, but you can
deduce it from the lengths of the edges on the boundary.

figure 2,12 shows the graph for the same member of
the second family.
\begin{center}
\resizebox{!}{1.4in}{\includegraphics{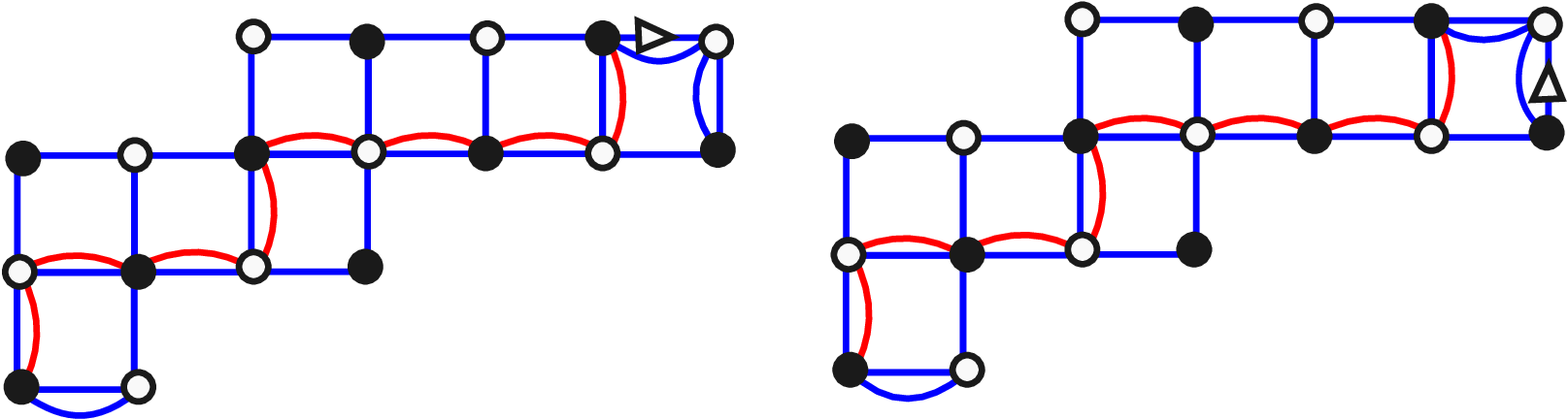}}
\newline
{\bf Figure 2,12:\/} Member $2324$ of the second staircase family
\end{center}

Figure 2.13 shows the corresponding nice coloring.

\begin{center}
\resizebox{!}{5in}{\includegraphics{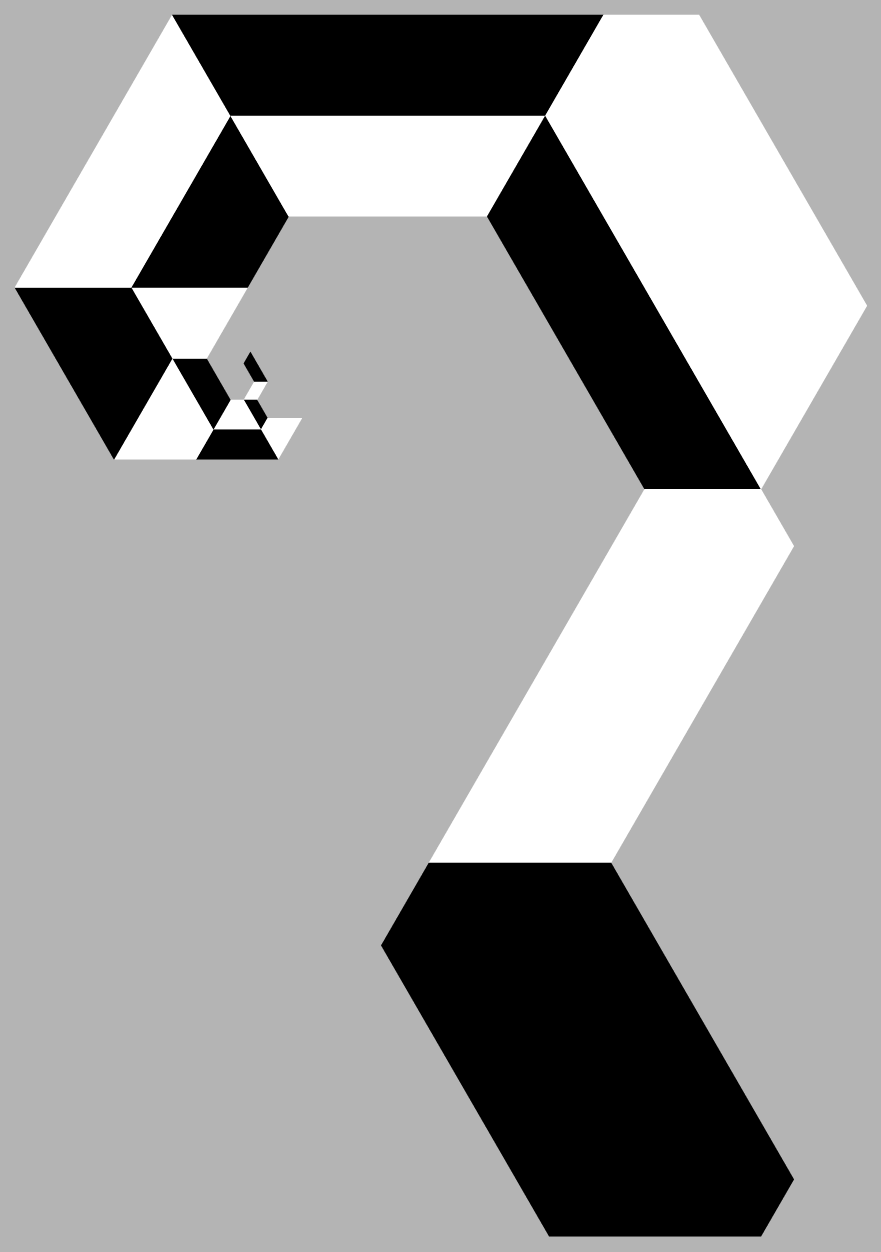}}
\newline
{\bf Figure 2.13:\/} The coloring associated to Figure 2.12.
\end{center}
Again, I am not showing the explicit folding pattern.
\newline
\newline
{\bf Remark:\/} In almost all the nice colorings I have
presented, the dual graph is made by gluing together
nearly identical disks.  There are plenty of examples
which do not have this feature, but I have not looked
for infinite families based on two very different halves.
The ``two halves'' approach is partly an artifact of
my program.  I have a feature where I make one
side, then duplicate it, then let the program find
the correct gluing of the disks to make a sphere.

\newpage

\section{Realizing Nice Colorings}
\label{realize}

\subsection{Axioms for the Multigraph}
\label{plausible}

Now we want to reverse the construction
in \S \ref{dualM}.
That is, we want
to start with an enhanced multigraph and
try to recover a partition of an octahedron
into nice polygons.
All our graphs will be planar multigraphs
which divide the sphere into quadrilaterals
and bigons.  Our multigraphs will have red and
blue edges.  We say that the {\it blue multigraph\/}
is the graph made entirely from the blue
edges.  We call a red edge and a blue edge
{\it parallel\/} if they share the same two
vertices.   Here are the conditions we need:
\begin{enumerate}
\item The blue multigraph should divide the
  sphere into $6$ bigons and some finite number
  of quadrilaterals.
\item Each quadrilateral face of the blue
  multigraph should contain one red edge
  that is parallel to one of the blue edges of the face.
\item The total number of edges incident to
  any vertex is $6$.
\end{enumerate}
We call a multigraph satisfying these conditions
a {\it plausible multigraph\/}.

We call the plausible multigraph {\it nice\/} if it
actually comes from a nice  coloring of an octahedron.
All the multigraphs
shown in the previous chapter are nice.
Given a plausible multigraph, we can extract the
shapes that would be involved in a nice coloring
of the sphere if such a thing existed.
Each vertex $v$ of blue degree $k$ in the plausible multigraph
corresponds to a nice polygon with $k$ sides. The
acute vertices correspond to faces in the (full) graph
which are incident to red edges.

\begin{center}
\resizebox{!}{2.1in}{\includegraphics{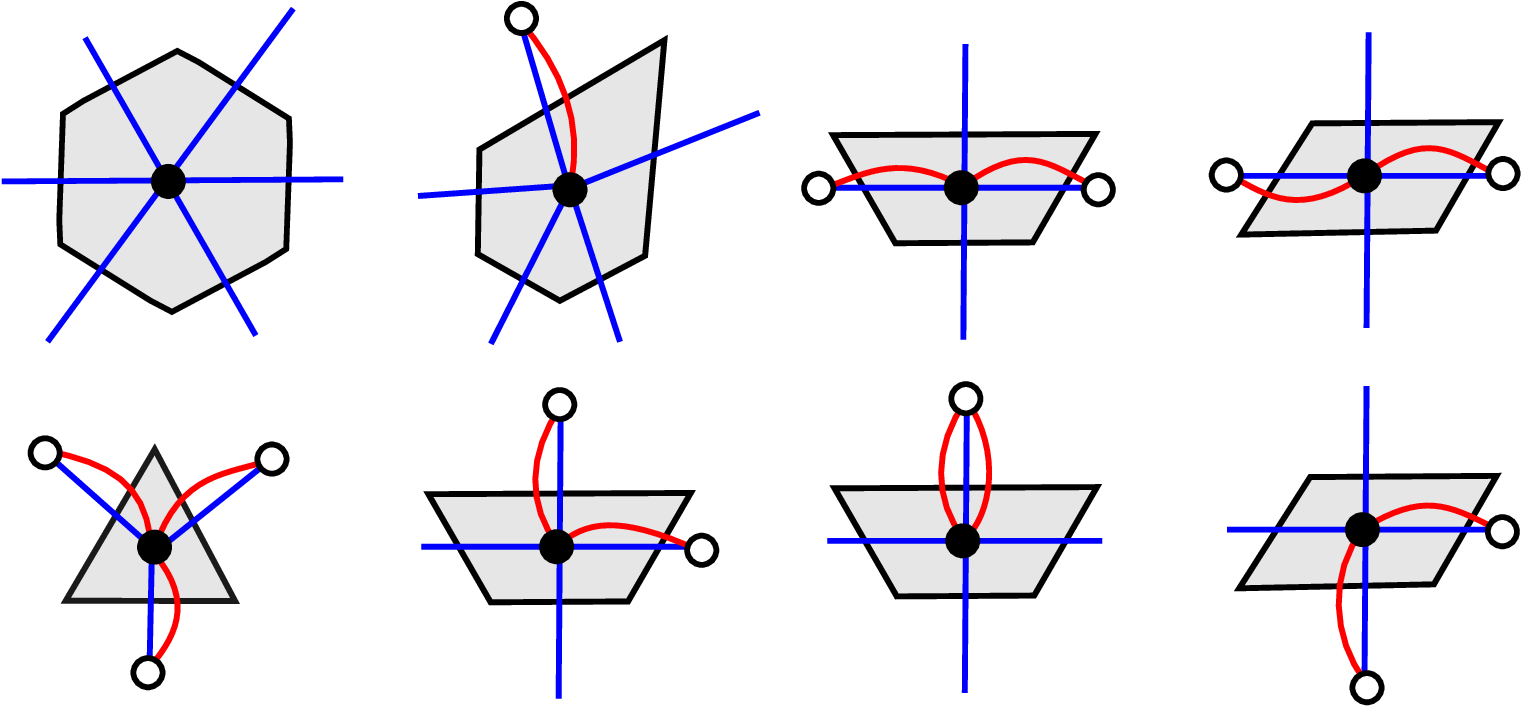}}
\newline
{\bf Figure 3.4:\/} Getting the shape of the nice polygon from the
graph
\end{center}

Figure 3.4 shows some
local pictures of the multigraph superimposed over
the corresponding shapes.  This is not quite an
exhaustive list of possibilities.

If $\widehat G$ is a nice multigraph, we can assign
a positive variable to each edge $e$ of $\widehat G$, namely
 the length of the polygon edge that $e$
crosses.  If $\widehat G$ is merely a plausible
multigraph, we want to find out if there exists
a positive assignment of variables to the edges
which corresponds to a nice coloring
of the sphere.    To do this, we solve a system of
linear equations and then see if the solution set
intersects the positive quadrant.  We now explain how
this is done.

\subsection{The Shape System}
\label{shapeS}

Let  $\omega=\exp(\pi i/3)$.
 We can view a
nice hexagon as a $6$-tuple of positive
numbers $\ell=(\ell_0,...,\ell_5)$ subject to the
following constraint.

\begin{equation}
  \label{closure2}
  \ell \cdot  (1,\omega,\omega^2,\omega^3,\omega^4,\omega^5)=0.
\end{equation}
Taking suitable linear combinations of the
real and imaginary parts we get the equivalent relations:
\begin{equation}
  \label{closure}
  \ell \cdot v_1=\ell \cdot v_2=0,
\hskip 10 pt
  v_1=(1,1,0,-1,-1,0), \hskip 10pt
  v_2=(0,1,1,0,-1,-1)
\end{equation}
The other nice polygons are limiting
cases of nice hexagons, where we set some of
the variables equal to $0$.  For instance,
if we have a trapezoid whose long side is
parallel to $\omega^0=1$ then we
set $\ell_1=\ell_5=0$ and then change the indices
so that they are consecutive. This gives
$\ell=(\ell_0,\ell_1,\ell_2,\ell_3)$ and 
$$
  \ell \cdot v_1=\ell \cdot v_2=0,
\hskip 10 pt
  v_1=(1,0,-1,0), \hskip 10pt
  v_2=(0,1,0,-1)
  $$

  \begin{lemma}
    \label{TRI}
  Suppose $P$ is a nice polygon having
  all integer sides.  Then $P$ has a triangulation
  by equilateral triangles with integer side
  lengths.
\end{lemma}

\startproof
When $P$ is a triangle the result is obvious.
Say that a {\it basic piece\/} is
a nice $4$-gon with integer side lengths.
If $P$ is a basic piece our proof goes
by induction on the sum of the side lengths of $P$.
In this case we can chop off an integer equilateral
triangle $T$ so that $P'=P-{\rm interior\/}(T)$
is a basic piece with smaller total side length.
The left side of Figure 3.5 shows this.

\begin{center}
\resizebox{!}{0.6in}{\includegraphics{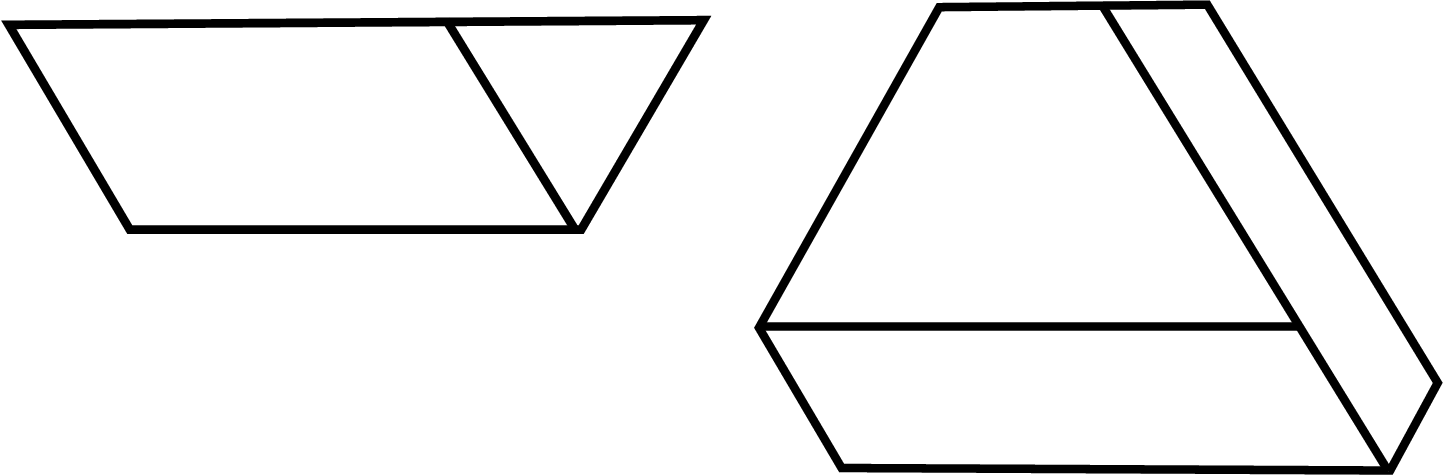}}
\newline
{\bf Figure 3.5:\/} An inductive proof
\end{center}

When $P$ is a pentagon,
we can chop off a basic piece and leave either
an integer equilateral triangle or a basic piece.  The idea here is to
use a basic piece which shares the shortest side of $P$.
When $P$ is a hexagon we can chop off a basic piece
and leave one of the pieces already considered,
thereby proceeding inductively.
The right hand side of Figure 3.5 shows how we reduce
a nice hexagon to a trapezoid with two chops.
\endproof

Let $\widehat G$ be a plausible multigraph, as above.
To (try to) find the nice colorings of the
sphere associated to $\widehat G$ we set up a system
of equations where we assign one real variable
to each blue edge of $\widehat G$.  The role played by
the red edges is that they help determine
the constraints.

We have $E_b$ variables, where $E_b$ is the
number of blue edges of $\widehat G$.
Each vertex of $\widehat G$ gives  two constraints
on the edges, as discussed in the previous section.
We represent the constraint as a vector of length
$E_b$ which is only nonzero in the positions
corresponding to the edges incident to the vertex.
We have $2V$ constraints, where $V$ is the
number of vertices of $\widehat G$.
We call this system of equations the {\it shape system\/}.
We call the shape system {\it nice\/} if it has
a positive solution.

Now let's compare the number $E_b$
of variables with the number $2V$ of
constraints.  We have
\begin{equation}
  \label{Euler}
  2F_b - 2E_b +2V = 4,
  \hskip 30 pt
  2F_b = E_b + 6, \hskip 20 pt
  \Longrightarrow \hskip 10 pt
  E_b-2V=2.
\end{equation}
The first equation is Euler's formula.   The last
equation, what we want, is a consequence of the first
two equations.  Here is how to derive the middle equation.
Set $F'_b=F_b-6$ and $E_b'=E_b-6$. These
respectively are the number of faces and
edges we get when we collapse each bigon
to a single edge.  That is, we identify the
edges bounding each bigon.   In the collapsed
graph we have $2F'_b=E'_b$, because each
face is a quadrilateral and each edge is incident
to $2$ faces.  Rewriting this in terms of
$F_b$ and $E_b$ we get the middle equation
in Equation \ref{Euler}.

\subsection{The Solution Space}
\label{solS}

Equation \ref{Euler} tells us that we have
two more variables than equations.
So, if we have one positive solution we
have (at least) a $2$-dimensional space
of positive solutions.
It turns out that the
constraint vectors are not independent.
In the next result we work with the
complex notation, as in Equation \ref{closure2}.
The next result, stated in the complex notation,
also says that there are $2$ real relations
amongst the constraint vectors.

\begin{lemma}
  \label{rel}
  The complex constraint vectors
  are linearly dependent.
  \end{lemma}

  \startproof
    The key to proving this result is to find
    a good basis for the constraint vectors.
    We will give the proof under the assumption
    that there exists a positive solution to the
    nice system.  The proof works in general,
    but in general it is harder to interpret the proof
    geometrically.

    Here is the key idea:  Let $\Sigma$ be the flat
    cone sphere which comes from a positive solution of
    the shape system.  There is a well defined
    {\it folding map\/} $f: \Sigma \to \C$.  The map
    is such that $f$ is an orientation preserving isometry
    on each white piece and an orientation reversing
    isometry on each black piece.

    Once we stipulate
    $f$ maps some particular edge into the real line,
    $f$ gives us a consistent way to label every edge
    in the dual graph $\widehat G$ by $6$th roots of unity.
    We orient the edges of the nice polygons in $\Sigma$
    so that they go clockwise around the white faces.
    Once we do this, we assign the $6$th root of unity
    $L(e)$ that is parallel to the image of the oriented edge $e$
    under $f$.  We  give this same label to the
    edge $e'$ of $\widehat G$ which crosses $e$.

    Now the constraint function for any given vertex $v$ in $\widehat G$ is
    \begin{equation}
      \label{constraint}
     \pm \sum_{j=1}^k L(e_j) \ell_j.
   \end{equation}
   Here
   \begin{itemize}
   \item $k$ is the number of blue edges incident to $v$,
   \item  $e_j$ is the $j$th such edge.
   \item $L(e_j)$ is the unit complex  label of $e_j$.
   \item $\ell_j$ is the length of the polygon edge crossed by $e_j$.
     \end{itemize}
       We choose $(+)$ if $v$ is a white vertex (corresponding to a white
    face)
    and $(-)$ if $v$ is a black vertex.

    Every edge label appears with a $(+)$ in one constraint vector
    and with a $(-)$ in another constraint vector and otherwise
    not at all.  Hence the sum of all the constraint vectors is $0$.
    This is our relation.
    \endproof

    \noindent
    {\bf Remark:\/} Even without a positive solution we
    could work out the labels of the edges of $\Sigma$.
    This is a purely algebraic construction, and the geometry
    only serves to interpret it.
    \newline

    The next lemma does not directly bear on
    Theorem \ref{one} because it discusses the complex rank
    of the set of constraint vectors rather than their real rank.
    For the examples computed in \S 2, the real rank is
    twice the complex rank.

    \begin{lemma}
      \label{RANK}
      The constraint vectors generate a
      $(V-1)$-dimensional complex subspace over $\C$.
    \end{lemma}

    \startproof
    Let $Y_1,...,Y_V$ denote the constraint vectors.
    We already know that $\sum Y_v$ is a relation.
    Suppose that $\sum c_v Y_v$ is another relation.
    We want to show that the constants
    $c_1,...,c_V$ are all the same.

    Each edge $e$ is incident to two vertices, say
    $v$ and $w$.  With Equation \ref{constraint} in
    mind, we see that the constraint vector $Y_x$
    has a $0$ in the $(e)$th coordinate unless
    $x \in \{v,w\}$.   The $(e)$th coefficients
    for $x=v$ and $x=w$ sum to $0$.  This means
    that $c_v=c_w$.  Since this is true for any edge
    $e$, we see that the coefficients of the relation
    are all equal. Hence this new relation is a multiple
    of the old one.
    \endproof

    \noindent
    {\bf Remark:\/} ChatGPT simplified the proof of
    Lemma \ref{RANK} for me. My original proof was
    somewhat less efficient.
    \newline

        Let me include the example shown to me
        by Kai Fu and Dia Taha. In this example, the
        two hexagons are glued around the boundary
        according to the pattern indicated by the arrows.

       \begin{center}
       \resizebox{!}{1.4in}{\includegraphics{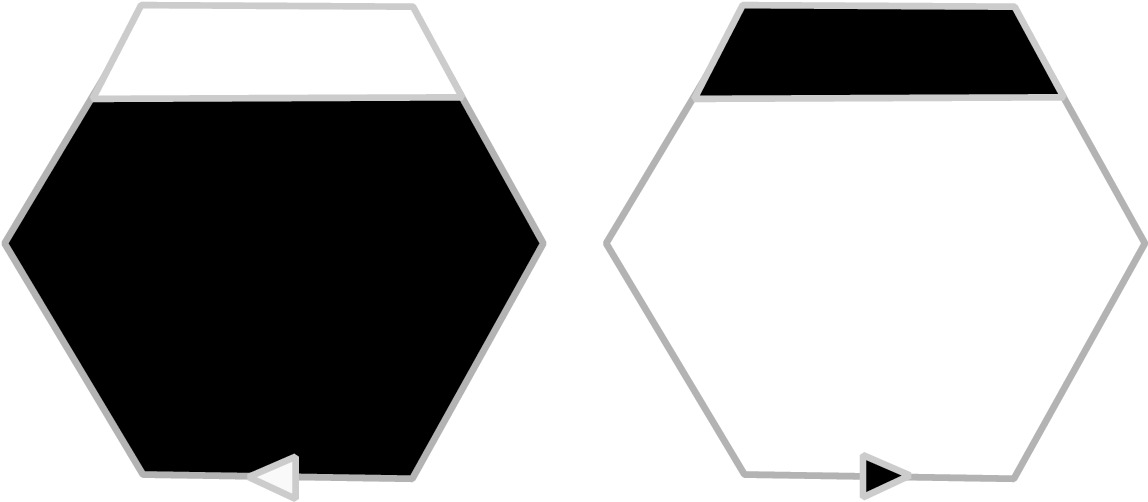}}
       \newline
       {\bf Figure 3.6:\/} An example with a $5$-dimensional moduli
       space
       \end{center}

        By adding horizontal
        cuts across the hexagons, and making the colors
        alternate on the two hexagons in a complementary way,
        you can increase the dimension as much as you like.

        \subsection{Putting it Together}
        \label{FINAL}

        Now we assemble the ingredients and prove
        Theorem \ref{one}.  Suppose that we have
        some nice coloring of the sphere.
        We let $\widehat G$ be the plausible multigraph that
        encodes the nice coloring.  The shape
        system has $E_b$ variables and, by Lemma
        \ref{rel}, has rank  at most $E_b-4$.   Hence the
        space of solutions to the shape system is at least
        $4$-dimensional.  Let $\cal S$ denote the
        solution space.

        The space $\cal S$ is a finite dimensional
        linear subspace of $\R^{E_b}$.  Also,
        $\cal S$ is defined by integer equations.
        So $\cal S$ is a rational subspace.
        We are
        interested in the subset where all variables
        are positive.  This is the intersection of
        $\cal S$ with the positive orthant, and
        this intersection is a rational convex
        polyhedral cone.  This is our cone $\cal C$
        from Theorem \ref{one}.

        Thanks to Lemma \ref{TRI}, each point of
        ${\cal C\/} \cap \Z^{E_b}$ corresponds to a
        triangulable nice coloring.  Again, since
        $\cal C$ is a rational cone, the linear
        subspace $\cal S$ intersects $\Z^{E_b}$ in a
        lattice.

        This completes the proof of Theorem \ref{one}.

\newpage

\section{The Map and the Quadratic Form}

In this chapter we prove
Theorem \ref{two}, which gives a locally affine
map from $\cal C$ to ${\cal M\/}$
and relates it to a quadratic form on $\cal C$.
We make the blanket assumption that we
are working with an example where
$\cal C$ is $4$-dimensional.

\subsection{Thurston Coordinates: A First Pass}
\label{triang}

We will first give an elementary account of
Thurston's coordinates on the space $\cal M$.
After this elementary
account we will give a second account that is
more conducive to the construction of our
linear map.  We give the first account mostly
for expository purposes.
\newline
\newline
\noindent
{\bf From Complex Numbers to Octahedra:\/}
As always, an {\it octahedron\/} for us is a
flat cone sphere with $6$ cone points having
cone angle $4 \pi/3$.
First we
describe how to build an octahedron
given complex numbers $z_0,z_1,z_2,z_3$.
Once $z_0$ is chosen, the other three
complex numbers need to be chosen with
some care, but an open set of choices works.
Let ${\bf Eis\/}$ denote the usual Eisenstein
lattice. This is the ring of integer combinations
of $1$ and $\omega=\exp(\pi i/3)$.
The number $z_0$ specifies a lattice of points
in the plane, namely $\Lambda=z_0 {\bf Eis\/}$.
What we are doing here is scaling up the
Eisenstein lattice by $z_0$.
We let $\Gamma$ be the group generated by the
order $3$ rotations about the elements of
$\Lambda$.  The quotient $\C/\Gamma$ is a
flat cone sphere with $3$ cone points with
cone angle $2\pi/3$.  This is a doubled
equilateral triangle.

Let $\Gamma_{\lambda}$ denote the order $6$ group
of rotations which fix $\lambda \in \Lambda$.
The stabilizer of $\lambda$ in $\Gamma$ has
index $2$ in $\Gamma_{\lambda}$.
Let $\lambda_1,\lambda_2,\lambda_3$
respectively be the points $0,z_0,\omega z_0$.
Notice that the three orbits
$\Gamma \lambda_k$ for $k=1,2,3$ are
pairwise disjoint, and their union is
all of $\Lambda$.

Around $\lambda_k$ we delete a hexagon
$h_k$ with the following description:
One of the vertices of $h_k$ is $\lambda_k + z_k$
and the others are the orbit of this one
vertex under $\Gamma_{\lambda_k}$.  We need to
choose $z_1,z_2,z_3$ so that the three hexagon
orbits $\Gamma h_k$ for $k=1,2,3$ are
pairwise disjoint, and indeed the grand union of
all the hexagons under the action of $\Gamma$
consists of disjoint hexagons. The sphere is
\begin{equation}
  \label{TC}
  \Sigma(z_0,z_1,z_2,z_3)=(\C-\Gamma h_1-\Gamma h_2-\Gamma h_3)/\Gamma.
\end{equation}

Why does this work?  As mentioned above,
the simpler quotient $\C/\Gamma$ is a flat cone
sphere with $3$ cone points with cone angle $2\pi/3$.
When we cut out each of the quotients
$\Gamma h_k/\Gamma$ we are cutting off
the cone point with cone angle $2\pi/3$ and replacing
it by $2$ cone points with cone angle
$4 \pi/3$.
\newline
\newline
{\bf Reversing the Process:\/}
Starting with an octahedron
$\Sigma \in {\cal M\/}$ we can group the $6$ cone points
in pairs and connect the two cone points in each pair
by an arc.  We can do the pairing and connecting so as
to minimize the total length of the three arcs.  Once we
do this, the arcs will be embedded.  (If the arcs cross, we
can surger them locally and find a new collection that is shorter.)

Let $\alpha_k$ for $k=1,2,3$ be the three arcs.
Since we have minimized the total length, each of
these arcs is a geodesic segment.
If we develop a small tubular
neighborhood of $\alpha_k$ into $\C$
we get a region of the form $U_k-h_k$ where $U_k$ is
some open topological disk and $h_k$ is
a hexagon with $6$-fold symmetry.
Figure 4.1 shows what we mean.

\begin{center}
\resizebox{!}{1.2in}{\includegraphics{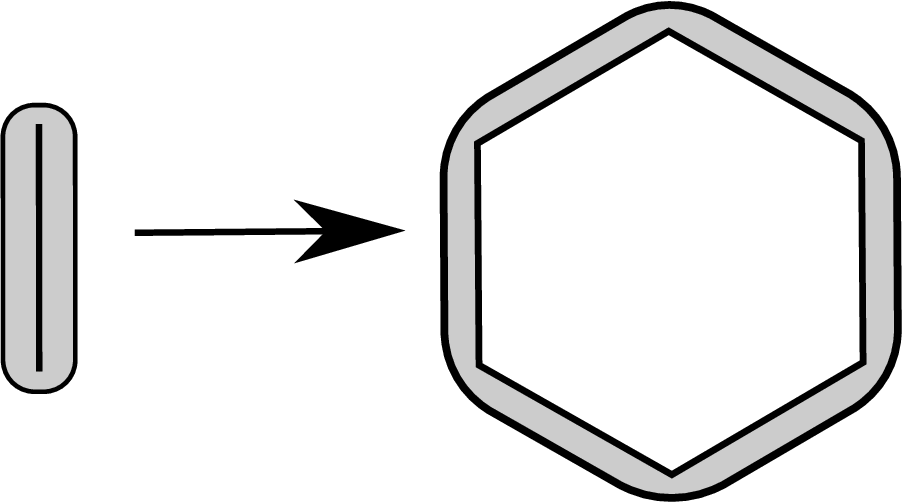}}
\newline
{\bf Figure 4.1:\/}  Developing the complement of an arc into $\C$
\end{center}

Let $\overline h_k$
be the quotient of $h_k$ by the order $3$ (and index $2$)
subgroup of
rotational symmetries of $h_k$.  We can isometrically glue
$\overline h_k$ to $\Sigma$ for $k=1,2,3$.
The result is a flat cone sphere $\Sigma'$ with $3$ cone
points all having cone angle $2\pi/3$.
The flat cone sphere $\Sigma'$ (after a choice) determines a lattice
$z_0{\bf Eis\/}$ and then the lifts of the glued
in copies of $\overline h_k$ determine the
hexagons $h_k$ for $k=1,2,3$ and the
complex numbers $z_1,z_2,z_3$.

There is one subtle point in this discussion. We say
{\it after a choice\/} because the geometry of
$\Sigma'$ only determines $z_0$ up to multiplication
by a unit complex number.  However, once we
make some choice, the collection of sufficiently
nearby $4$-tuples $(z_0,z_1,z_2,z_3)$ will parametrize
a neighborhood of $\Sigma$ in ${\cal M\/}$.
In this moduli space, the spheres
$\Sigma(z_0,z_1,z_2,z_3)$ and
$\Sigma(uz_0,uz_1,uz_2,uz_3)$ are considered distinct
even though, when $u$ is a unit complex number, the
underlying spheres are isometric.

These coordinates just give a
local coordinate system, for several reasons. First, we had to make
a choice of how to pair up the cone points based on the
geometry of $\Sigma$.   Second, if we replace $z_k$ by $\omega z_k$ we
get a different parametrization of the same sphere. So, we
can only change the parameters a little bit if we want
bijective parametrization of a neighborhood of $\Sigma$
in ${\cal M\/}$.

\subsection{Thurston Coordinates: A Second Pass}
\label{coord2}

Again following [{\bf T\/}], 
we take a more sophisticated approach to Thurston's
coordinates on $\cal M$. The new system we get
is equivalent to the old one {\it via\/} a complex
linear change of coordinates.

Let $\Sigma$ be a flat cone sphere.
Let $c$ denote the set of $6$-cone
points of $\Sigma$ and let
$\Sigma'=\Sigma-c$.    Let
$\widetilde \Sigma'$ denote the universal
cover of $\Sigma'$.  Let
$\pi: \widetilde \Sigma' \to \Sigma'$ denote
the universal covering map.
The space $\widetilde \Sigma'$ is not metrically
complete.  We can have a Cauchy sequence
$\{\widetilde p_n\}$ of points
in $\widetilde \Sigma'$ such that
$\{p_n\}$ converges to a cone point.
Here we have set $p_n=\pi(\widetilde p_n)$.
The limit of $\{\widetilde p_n\}$ will not
exist in $\widetilde \Sigma'$.

We let $\widetilde \Sigma$ denote the metric
completion of
$\widetilde \Sigma'$.
The space $\widetilde \Sigma$
is known as the
{\it orbifold universal cover\/} of $\Sigma$.
We get $\widetilde \Sigma$ from $\widetilde \Sigma'$
by adding in points which serve
as the limits of all the Cauchy sequences
considered in the previous paragraph.
The difference
$\widetilde \Sigma-\widetilde \Sigma'$ is a
discrete countable set of points which
$\pi$ maps to $c$.  We call these points
the {\it special points\/}.

The following is a lemma from [{\bf T\/}].

\begin{lemma}
  $\widetilde \Sigma$ has an equivariant triangulation whose
  edges are geodesic segments having special points as
  endpoints.
\end{lemma}

\startproof
 Here is a sketch. Each special
point $p$ defines a {\t Voronoi cell\/} $K_p$ consisting
of the closure of all the points in $\widetilde \Sigma$
which are closer to $p$ than to any other special point.
About
each point in the intersection of at least $3$ Voronoi
cells is a metric disk that contains a finite number
of points in its boundary.  The convex hull of these
points is a convex polygon whose vertices are
special points.  These polygons are
called {\it Delaunay cells\/}.  Typically the Delaunay
cells are triangles; in case they are not, they may
be further triangulated into triangles by the addition
of some of the diagonals. This further triangulation
may be done in an equivariant way.
\endproof

Since the triangulation of $\widetilde \Sigma$ is equivariant,
it induces a triangulation of the original flat cone
sphere $\Sigma$.  We can choose a collection
of $4$ edges $e_1,e_2,e_3,e_4$ of the triangulation of $\Sigma$ which
make a tree, and then lift this tree to $\widetilde \Sigma$.
To be artistic about it, we can choose the tree to look
like a quadrapod, with the central vertex being equal to $v$.
Then we can take the lifted tree to have $\widetilde v$ as
the corresponding lifted vertex.

Let $\widetilde w_k$ be the other endpoint of the
lift $\widetilde e_k$ of $e_k$. The $4$ complex
numbers $d(\widetilde w_k)$ for $k=1,2,3,4$ serve as local
coordinates for the space $\cal M$.

\subsection{The Affine Map}

Now we describe a real affine map from
$\cal C$ into $\C^4$ and then we interpret
this map as a locally affine map from $\cal C$ into
$\cal M$.

If we have a nice coloring of $\Sigma$, we can lift
this to get a nice coloring of $\widetilde \Sigma$.
To normalize our constructions, we fix a triple
$(v,e,f)$ where $f$ is a polygon
of the nice coloring, and $e$ is an edge of $f$, and
$v$ is a vertex of $f$ incident to $e$.  We also
suppose that $v$ is a cone point.
Now we choose a lift
$(\widetilde v,\widetilde e,\widetilde f)$ of
$(v,e,f)$ to $\widetilde \Sigma$.  We call this
lift our {\it favorite flag\/}.
We have a unique locally isometric map
$d: \widetilde \Sigma \to \C$ such that
\begin{itemize}
\item $d(\widetilde v)=0$.
\item $d(e) \subset \R$ and $d(e-v) \subset \R^+$.
\item $d(f-e) \subset \R \times \R_+$.
\end{itemize}
Figure 4.2 shows our normalization relative to
the coordinate axes in $\C$.
  
\begin{center}
\resizebox{!}{1.5in}{\includegraphics{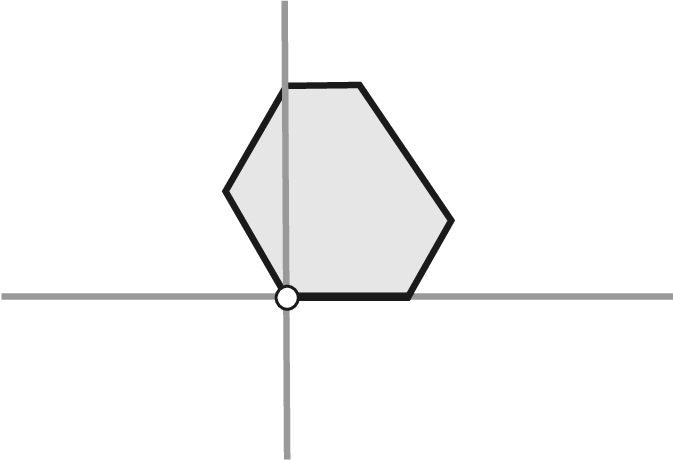}}
\newline
{\bf Figure 4.2:\/}  Normalizing the developing map.
\end{center}

We pick a point in
$\cal C$.  This defines a nice coloring for us.
We then build $\widetilde \Sigma$ based on
this coloring, normalize as above, and
then take the point in
$\cal M$ with coordinates
$d(w_k)$ for $k=1,2,3,4$.
This is our map $A$.

The map from $\cal C$ to $\C^4$ is real linear.
To see this, note that there is some finite
chain of nice polygons $P_1,...,P_m$ which
cover the geodesic segment connecting
our favorite vertex $\widetilde v$ to the
vertex $\widetilde w_k$.   We can express
$d(w_k)$ as the sum of certain complex
numbers describing certain of the edges of
$d(P_j)$ for $j=1,...,m$.  When we
change the edge labels for our coloring,
this sum changes in a linear way.

One thing we note is that the map
  $A({\cal S\/}) \subset \C^4$ to
  $\cal M$ need not be injective.
  Since we are working in a local coordinate
  chart, all we can say is that $A: {\cal C\/} \to {\cal M\/}$
  is locally affine and locally injective.

\subsection{The Thurston Hermitian Form}

Thurston introduces a Hermitian form on the space
$\cal M$.
In the local coordinates described in
\S \ref{triang}, the form is given by

\begin{equation}
  \label{Hform}
  \langle Z,W\rangle = 6z_0 \overline w_0 -  6z_1 \overline w_1 - 6z_2
  \overline w_2 -  6z_3 \overline w_3.
\end{equation}
Here $Z=(z_0,z_1,z_2,z_3)$ and $W=(w_0,w_1,w_2,w_3)$.
The diagonal part of this form $Z \to \langle Z,Z\rangle$
computes a multiple of the area of $\Sigma$.
The multiple is such that the form counts $3$ times
the number of triangles when the tiling is by
equilateral triangles having side length $1$.

In general, the changes of coordinates between Thurston's
various coordinate systems are always complex linear.
The Hermitian form has a local expression in each coordinate
system.   The secret to lining up all the different local
expressions is that their diagonal parts all compute the
area, and a Hermitian form on a complex vector space
is determined by its diagonal part.
Equation \ref{Hform} 
reveals that the form has signature $(1,3)$.

Let me describe what the
Hermitian form looks like in the coordinate
system defined in \S \ref{coord2}.
This alternate perspective connects with the
way we define the quadratic form on $\cal C$
in the next section.
We have triangulated
the universal cover $\widetilde \Sigma$.
Each triangle in the triangulation belongs
to an orbit.  We choose a finite collection
of triangles, one per triangle orbit.  We
then develop these triangles in the plane.
The locations of the points will be complex
linear functions of $d(e_k)$.

Suppose we have two such triangulations,
corresponding to different points of
$\cal M$.  Let $\tau$ and $\tau'$ be
two corresponding triangles.  We choose
a vertex $v$ of $\tau$ and the corresponding
vertex $v'$ of $\tau'$.  We let
$z_1$ and $z_2$ be the complex numbers
which describe the difference between the
other vertices of $\tau$ and $v$.
Likewise define $z_1'$ and $z_2'$.  For
ease of notation we set $w_k=z_k'$.
We then consider the expression
\begin{equation}
  \label{localform}
  \frac{1}{4i}\bigg( z_1\overline{w_2} - z_2 \overline{w_1}\bigg)
\end{equation}
When $\tau=\tau'$ this expression computes the
area of $\tau$.
In these coordinates, the Hermitian form is the
sum of these expressions, appropriately scaled
to match Equation \ref{Hform}.

\subsection{The Quadratic Form}
\label{QF}

Suppose we have a nice hexagon with
consecutive side lengths $\ell=(\ell_1,...,\ell_6)$.  Then,
\begin{equation}
  Q(\ell,\ell) = 2\sum_i\ell_i \ell_{i+1}+
  \sum_i \ell_i \ell_{i+2}
\end{equation}
computes $8 \sqrt 3$ times the area.
The indices are taken cyclically.
You can derive this by the following procedure:
\begin{itemize}
\item Triangulate the hexagon.
\item Compute the areas of the triangles using Equation
  \ref{localform}
\item Symmetrize the resulting expression by averaging over all cyclic
  permutations of the variables.
\end{itemize}
We get similar formulas for other nice
polygons by setting various of the variables
equal to $0$.

With the understanding that we are writing just
one formula for all combinatorial types of nice
polygon, we introduce the real quadratic form.
\begin{equation}
  \label{niceform}
  Q(\ell,m) = 2 \sum_i(\ell_i m_{i+1}+ \ell_{i+1} m_i) +
 \sum_i(\ell_i m_{i+2} +  m_i \ell_{i+2})
\end{equation}
We will only apply this form for nice polygons
having the same combinatorics.
To explain the normalization, we note that when this
formula is applied to the pair $\ell=m=(1,1,1,1,1,1)$
it returns $18$, which is $3$ times
he number of unit equilateral triangles
tiling the corresponding unit hexagon.
The reason for the normalization is that we wanted
to get an integral quadratic form.

Supposing that our nice coloring has $n$ nice polygons,
we define the following quadratic form on
$\cal C$:
\begin{equation}
  \label{sum}
  Q = \sum_{i=1}^n Q_i,
\end{equation}
where $Q_i$ is the quadratic form associated
to the $i$th nice polygon.
This is our quadratic form on $\cal C$.
This normalization computes $3$ times
the number of
triangles in $\Sigma$ when $\Sigma$ is triangulable
by equilateral triangles of side length $1$.

\subsection{Non-Degeneracy and Injectivity}

We have our real linear map $A: {\cal C\/} \to \C^{1,3}$.
This map extends to a real linear map
$A: {\cal S\/} \to L \subset \C^{1,3}$.
Here $\cal S$ is the $4$-dimensional
subspace spanned by $\cal C$ and
$L$ is a $4$-dimensional real linear subspace.
At least on an open neighborhood of the
good coloring on which this whole apparatus is
based, the map from $L$ into $\cal M$ is
injective.

Let us define $Q'$ to be $6$ times the real part of
Thurston's Hermitian form, restricted to $L$.
The diagonal part of $Q'$ computes $(8 \sqrt 3)$ times the area,
at least when restricted to subsets of $L$
corresponding to points in $\cal M$.
This is the same factor that comes up for $Q$.

Let $A^*$ be the pullback operator for
quadratic forms.  The two forms
$A^*(Q')$ and $Q$  have the
same diagonal parts: on suitable open
sets, both compute
area.  But a quadratic form on a real
vector space is determined by its
diagonal part. Hence $A^*(Q')=Q$.

\begin{lemma}
  If $Q$ is non-degenerate then
  $A$ is locally injective.
\end{lemma}

\startproof
Suppose for the sake of contradiction
that $A: {\cal S\/} \to L$
has a kernel.  Then for any nonzero
vector $W \in \cal S$ in the kernel of $A$ we
have $A^*(Q')(W,W)=0$.  This contradicts
the fact that $A^*(Q')=Q$ and that $Q$ is
non-degenerate.
\endproof

\begin{lemma}
  Suppose that $Q$ has signature $(1,3)$.
  Then the projectivized map $PA$ is
  locally injective on $P\cal C$.
\end{lemma}

\startproof
If $PA$ is not locally injective,
then there are two points
$V,W \in L$ such that
$W=\lambda V$ where $\lambda \in \C$
is non-real.  But then
$L$ contains the $\C$-span of $V$.
We have
$$Q'(V,V)>0, \hskip 30 pt Q'(iV,iV)>0, \hskip 30 pt
Q'(V,iV)=0.$$
The then the restriction of $Q'$ to the $\C$-span of
$V$ has signature $(2,0)$.    This means that
$Q'$ cannot have signature $(1,3)$ on $L$.
\endproof

\newpage

.

\section{References}

\noindent
[{\bf BG\/}] C. Bavard and E. Ghys, {\it Polygones du plan et
  poly\`edres hyperboliques\/}, Geom. Dedicata 43(2),  (1992) pp 207--224
\vskip 9 pt
\noindent
[{\bf T\/}] W. Thurston, {\it Shapes of Polyhedra and triangulations
  of the sphere\/}, \newline
Geometry $\&$ Topology Monographs, Vol. 1:
The Epstein birthday shrift, pp 511-549  (1998)  See also
arXiv:math/9801088.
\vskip 9 pt
\noindent
[{\bf W\/}] D. B. West, {\it Introduction to Graph Theory, 2nd
  Edition\/}, Prentice Hall (2001)

\end{document}